\documentclass[11pt, letterpaper, oneside]{compositio}

\usepackage{latexsym, amsmath, amssymb, amsfonts, amscd}
\usepackage{amsthm}
\usepackage{t1enc}
\usepackage[mathscr]{eucal}
\usepackage{indentfirst}
\usepackage{graphicx, pb-diagram}
\usepackage{fancyhdr}
\usepackage{fancybox}
\usepackage{enumerate}
\usepackage[all]{xy}

\theoremstyle{plain}
\newtheorem{thm}{Theorem}[section]

\newtheorem{prop}[thm]{Proposition}
\newtheorem{lemma}[thm]{Lemma}

\theoremstyle{definition}
\newtheorem{defi}[thm]{Definition}
\newtheorem{pdef}[thm]{Proposition-Definition}

\theoremstyle{remark}
\newtheorem{remark}{Remark}
\newtheorem{ep}{Example}

\newtheorem{note}{Note}

\newcommand{\lra}{\longrightarrow}

\newcommand{\btd}{\bigtriangledown}

\newcommand{\Z}{\ensuremath{\mathbb Z}}

\newcommand{\R}{\ensuremath{\mathbb R}}

\newcommand{\di}{d}
\newcommand{\ci}{\ensuremath{C^{\infty}}}

\newcommand{\cC}{\mathcal{C}}            
\newcommand{\cF}{\mathcal{F}}
\newcommand{\cX}{\mathcal{X}}
\newcommand{\cY}{\mathcal{Y}}
\newcommand{\cZ}{\mathcal{Z}}
\newcommand{\cG}{\mathcal{G}}
\newcommand{\cH}{\mathcal{H}}

\newcommand{\half}{\frac{1}{2}}

\newcommand{\bt}{\mathbf{t}}                  
\newcommand{\bs}{\mathbf{s}}                  
\newcommand{\bbt}{\bar{\mathbf{t}}}           
\newcommand{\bbs}{\bar{\mathbf{s}}}           
\newcommand{\bm}{\bar{m}}                     

\begin{document}

\title{Integrating Lie algebroids via stacks}
\author{Hsian-Hua Tseng}
\email{hhtseng@math.berkeley.edu}
\address{Department of Mathematics, University of California, Berkeley, CA
94720}

\author{Chenchang Zhu}
\email{zcc@math.berkeley.edu}
\address{Department of Mathematics, University of California, Berkeley, CA
94720} 

\date{\today}

\begin{abstract}
Lie algebroids can not always be integrated into Lie groupoids. We
introduce a new object---``Weinstein groupoid'', which is a
differentiable stack with groupoid-like axioms. With it, we present a
soloution to the integration problem of Lie algebroids. It turns out
that every Weinstein groupoid has a Lie algebroid and every Lie
algebroid can be integrated into a Weinstein groupoid.
\end{abstract}

\maketitle

\section{Introduction}
In this paper, we present a new viewpoint to integrate Lie algebroids:
unlike (finite dimensional) Lie algebras\footnote{Non-integrability
already appears in the case of infinite dimensional Lie algebras
\cite{dl}. In this paper, Lie algebroids are assumed to be finite
dimensional. }  which always have their associated Lie groups, Lie
algebroids do not always have their associated Lie groupoids
\cite{am1} \cite{am2}. So the Lie algebroid version of Lie's third
theorem poses the question indicated by the following chart:

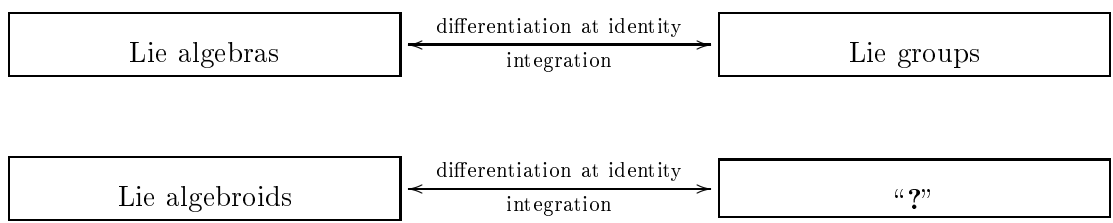
\[ \xymatrix{
& \fbox{\parbox{.3\linewidth}{\center{Lie algebras}}}
  \ar[rrrr]^{\text{differentiation at identity}} &  & & &
  \fbox{\parbox{.3\linewidth}{\center{Lie groups}}}
  \ar[llll]^{\text{integration}}   \\
 & \fbox{\parbox{.3\linewidth}{\center{Lie
   algebroids}}} \ar[rrrr]^{\text{differentiation at identity}} &  & & &
   \fbox{\parbox{.3\linewidth}{\center{``{\bf ?}''}}}
  \ar[llll]^{\text{integration}} }
\]

Pradines posed the above question in \cite{pradines} and constructed
local Lie groupoids (formulated in \cite{cdw} \cite{karasev}
\cite{van-est}) as the integration object ``{\bf ?}''. But a global object for ``{\bf ?}'' is
still in need: not only it gives conceptually better answer to the
diagram above (Lie groups are globle objects), but also it has
profound applications in Poisson geometry, such as Weinstein's
symplectic groupoids \cite{w}, Xu's Morita equivalence of Poisson
manifolds \cite{x1} \cite{x2}, symplectic realizations \cite{w-poisson},
Picard groups \cite{picard} and the linearization problem of Poisson manifolds
\cite{cf3}.

After Pradines' local groupoids, progress towards special cases of the
above integration problem was made, among others by \cite{dazord1}
\cite{debord} \cite{Mac}
\cite{nistor} \cite{w-heis}. In \cite{cafe}, Cattaneo and Felder realized the symplectic groupoid
as the phase space of the Poisson sigma model, where path spaces were
natually used. In the recent work \cite{cf}, Crainic and Fernandes
established the $A$-path space, which has become a widely used method
in groupoid theory. With the $A$-path space, they constructed a
topological groupoid for every Lie algebroids and finalized the
integrability condition for Lie algebroids. Weinstein further
conjectured that this topological groupoid has some differentiable
structure \cite{picard}. But ordinary structures such as manifolds and
orbifolds don't serve the purpose. It turns out that differentiable
stacks discussed in recent papers \cite{bx} \cite{metzler}
\cite{pronk} provide a suitable structure to above conjecture posed by
Weinstein.

Based on all this previous work, we further extend the research by
defining a more general concept---Weinstein groupoids. Roughly
speaking, they are groupoids in the world of stacks (see Definition
\ref{wgpd}). Weinstein groupoids are global differentiable objects
providing a soulution to the Lie algebroid version of Lie's third
theorem.

\begin{defi} [Weinstein groupoid]\label{wgpd}
A Weinstein groupoid over a manifold $M$ consists of the following
data:
\begin{enumerate}
\item an \'etale differentiable stack $\cG$ (see Definition \ref{stack});

\item (source and target) maps $\bar{\bs}$,
$\bar{\bt}$: $\cG \to M$ which are surjective submersions between
differentiable stacks;

\item (multiplication) a map $\bar{m}$: $\cG\times_{\bbs, \bbt} \cG \to
\cG$, satisfying the following properties:
\begin{itemize}
  \item $\bbt \circ \bar{m}=\bbt\circ pr_1$, $\bbs \circ \bar{m}=\bbs\circ
  pr_2$, where $pr_i: \cG \times_{\bbs, \bbt} \cG \to \cG$ is the
  $i$-th projection $\cG\times_{\bbs, \bbt} \cG \to
\cG$;
  \item associativity up
to a 2-morphism, i.e. there is a unique 2-morphism $\alpha$
between maps $\bar{m}\circ (\bar{m} \times id)$ and
$\bar{m}\circ(id\times \bar{m})$;
\end{itemize}

\item  (identity section) an injective immersion  $\bar{e}$: $M\to \cG$
such that, up to 2-morphisms, the following identities
\[
\bar{m}\circ ((\bar{e}\circ \bbt)\times id)=id, \,\,\bar{m}\circ
(id\times (\bar{e}\circ\bbs) )=id,\] hold (In particular,by
combining with the surjectivity of $\bbs$ and $\bbt$, one has
$\bbs \circ \bar{e}= id$, $\bbt \circ \bar{e}= id$ on $M$);

\item (inverse) an isomorphism of differentiable stacks
$\bar{i}$: $\cG \to \cG$ such that, up to 2-morphisms, the following
identities
\[ \bar{m}\circ (\bar{i}\times id \circ \Delta)=\bar{e}\circ\bbs, \;\;
\bar{m}\circ (id\times\bar{i}\circ \Delta)=\bar{e}\circ \bbt,\]
hold, where $\Delta$ is the diagonal map: $\cG\to \cG\times\cG$.
\end{enumerate}
Moreover, restricting to the identity section, the above
2-morphisms between maps are the $id$ 2-morphisms. Namely, for
example, the 2-morphism $\alpha$ induces the $id$ 2-morphism
between the following two maps:\[ \bar{m}\circ ((\bar{m} \circ
(\bar{e}\times\bar{e}\circ \delta))\times \bar{e} \circ
\delta)=\bar{m}\circ(\bar{e}\times(\bar{m}\circ(\bar{e}\times\bar{e}\circ\delta))\circ\delta),
\]where $\delta$ is the diagonal map: $M\to M\times M$.
\end{defi}

\noindent{\em General Remark}: the terminology involving stacks in
the above definition, as well as in the following theorems, will
be explained in detail in Chapter \ref{sect: stack}. For now,  to
get a general idea of these statements, one can take stacks simply
to be manifolds.

Our main result is the following theorem:

\begin{thm}[Lie's third theorem] \label{lieIII}
To each Weinstein groupoid one can associate a Lie algebroid. For
every Lie algebroid $A$, there are naturally two Weinstein
groupoids $\cG(A)$ and $\cH(A)$ with Lie algebroid $A$.
\end{thm}

We can apply our result to the classical integrability problem,
which studies when exactly a Lie algebroid can be integrated into
a Lie groupoid.

\begin{thm}\label{integ}
A Lie algebroid $A$ is integrable in the classical sense iff
$\cH(A)$ is representable, i.e. it is an honest (smooth) manifold. In this
case $\cH(A)$ is the source-simply connected Lie groupoid of $A$
(it is also called the Weinstein groupoid of $A$ in \cite{cf}).
\end{thm}

We can also relate our work to the previous work on the
integration of Lie algebroids via the following two theorems:

\begin{thm}\label{local}
Given a Weinstein groupoid $\cG$, there is an\footnote{It is canonical up to isomorphisms near the
identity section.} associated local Lie groupoid $G_{loc}$ which
has the same Lie algebroid as $\cG$.
\end{thm}

\begin{thm}\label{cf}
The orbit spaces of $\cH(A)$ and $\cG(A)$ as topological spaces
are both isomorphic to the universal topological groupoid of $A$
constructed in \cite{cf}.
\end{thm}

\section{Path spaces}

We define the $A_0$-path space, which is very
similar to\footnote{Actually it is a submanifold of the $A$-path
space.} the $A$-paths defined in \cite{cf}.

\begin{defi} \label{apath}
Given a Lie algebroid\footnote{Here we require $A$ to be a Hausdorff
  manifold. } $A \overset{\pi} {\lra} M$ with  anchor
$\rho : A \to TM$, a $C^1$ map $a$: $I=[0, 1] \to A$ is an {\bf
$A_0$-path} if
\[ \rho (a(t)) = \frac{\di}{\di t} \left( \pi \circ a(t) \right) , \]
with the following boundary conditions,
\[ a(0)=0,\; a(1)=0,\; \dot{a}(0)=0,\; \dot{a}(1)=0.\]
We often denote the base path $\pi\circ a(t)$ in $M$ by
$\gamma(t)$. We denote $P_0A$ the set of all $A_0$-paths of $A$.
It is a topological space with topology given by uniform
convergence of maps. Omitting the boundary condition above, one
get the definition of {\bf $A$-paths}, and we denote the space of
$A$-paths by $P_aA$.
\end{defi}

We can equip $P_0A$ with the structure of a smooth (Banach)
manifold using a Riemannian structure on $A$. On the total space
of $C^1$ paths $PA= C^1(I, A )$, there is a $\ci$-structure as
follows: at every point $a: I \to A$ in $PA$, let $a^*TA \to I$ be
the pull-back of the tangent bundle to $I$. For $\epsilon>0$, let
$T_{\epsilon}\subset a^*TA$ be the open set consisting of tangent
vectors of length less than $\epsilon$. For sufficiently small
$\epsilon$, we have the exponential map $\exp$: $T_{\epsilon} \to
A$, $(t, v) \mapsto \exp_{a(t)} v$. It maps $T_{\epsilon}$ to an
open subset of $A$. Using this map we can identify $PT_{\epsilon}$
,the  $C^1$-sections of $T_{\epsilon}$, with an open subset of
$PA$. The oriented vector bundle $a^*TA$ over $I$ is trivial. Let
$\varphi : a^*TA \to I\times\R^n$ be a trivialization where $n$ is
the dimension of $A$. Then $\varphi$ induces a mapping from
$PT_{\epsilon}$ to $P\R^n=C^1(I, \R^n)$. Since $C^1(I, \R^n)$ is a
Banach space with norm $\|f\|^2=\sup\{|f|^2+ |f'|^2\}$,
$PT_{\epsilon}$ can be used as a typical Banach chart for the
Banach manifold structure of $PA$. $P_0A$ is defined by equations
on $PA$ which, in above local charts, can be written as
$$ \dot{\gamma}^k(t)= \sum_{j=1}^{m-n} \rho^k_j (\gamma(t)) a^j(t),
\;\;\; a^j(0)=a^j(1)=0, \; \dot{a}^j(0)=\dot{a}^j(1)=0, $$ for
$j=1, ..., n=\text{rank} A$, $k=1,..., m=\dim M$. 
The space of
$\{a^j\}$ satifying the boundary conditions is  a closed subspace
$C^1_0(I,\R^k)$  of $C^1(I, \R^k)$,  hence is also a
Banach space. Then $\gamma^i(t)$ is determined
by the initial value $\gamma^i(0)$ and $\dot{\gamma}^i(0)=0$. Hence
$P_0 A$ is a Banach manifold with a typical local chart the Banach
space $C^1_0(I,\R^k) \times \R^{n-k}$. We
refer to \cite{lang} for the definition and further properties of
Banach manifolds.


.

\begin{pdef}\label{homotopy} Let $a(\epsilon, t)$ be a family of
$A_0$-paths of class $C^2$ in $\epsilon$ and assume that their
base paths $\gamma(\epsilon,t)$ have fixed end points. Let
$\nabla$ be a connection on $A$ with torsion $T_{\nabla}$ defined
as
\[ T_{\nabla} ( \alpha , \beta) = \nabla_{\rho(\beta)} \alpha -
\nabla_{\rho(\alpha)} \beta + [\alpha, \beta]. \] Then the
solution $b=b(\epsilon, t)$ of the differential
equation\footnote{Here, $T_{\nabla} (a, b)$ is not quite well
defined. We need to extend $a$ and $b$ by sections of $A$,
$\alpha$ and $\beta$, such that $a(t)=\alpha(\gamma(t), t)$ and
the same for $b$. Then $T_{\nabla} (a, b)|_{\gamma(t)}:=$ $
T_{\nabla} (\alpha, \beta)|_{\gamma(t)}$ at every point on the
base path. However, the choice of extending sections does not
affect the result.}
\begin{equation}\label{eq-homotopy}
\partial_t b -\partial_{\epsilon }a = T_{\nabla} (a, b),
\;\;\;\;\; b(\epsilon ,0)=0 \end{equation} does not depend on the
choice of connection $\nabla$. Furthermore, $b(\cdot, t)$ is an
$A$-path for every fixed $t$, i.e. $\rho (b(\epsilon, t)) =
\frac{\di}{\di \epsilon} \gamma(\epsilon, t)$. If the solution $b$
satisfies $b(\epsilon,1)=0$, for all $\epsilon$, then $a_0$ and
$a_1$ are said to be  {\bf equivalent} and we write $a_0 \sim
a_1$.
\end{pdef}
\begin{remark} A homotopy of $A$-paths \cite{cf}  is
  defined by replacing $A_0$ by $A$ in the definition above. A similar result as above holds for $A$-paths \cite{cf}.  So the above statement holds viewing $A_0$-paths as $A$-paths.
\end{remark}

This flow of $A_0$-paths $a(\epsilon, t)$ generates a foliation
$\cF$.  The $A_0$-path
space is a Banach submanifold of the $A$-path space and $\cF$ is the
restricted foliation of the foliation defined in Section 4 of
\cite{cf}. For any foliation, there is an associated {\bf monodromy groupoid} \cite{moerdijk} (or {\bf fundamental
groupoid} as in \cite{cw}) : the objects are points in the
manifold and the arrows are paths within a leaf up to homotopies
with fixed end points inside the leaf. The source and target maps
associate the equivalent class of paths to the starting and ending
points respectively. For any regular foliation on a smooth
manifold its monodromy groupoid is a Lie groupoid in the sense of
\cite{cf}. In our case, it is an infinite dimensional groupoid
equipped with a Banach manifold structure. Here, we slightly
generalized the definition of Lie groupoids to the category of
Banach manifolds by requiring the same conditions but in the
sense of Banach manifolds. Denote the monodromy groupoid of $\cF$
by $Mon(P_0 A) \underset{\mathbf{t}_M}{ \overset{\mathbf{s}_M}
{\rightrightarrows}} P_0 A$. In a very similar way
\cite{moerdijk}, one can also define the {\bf holonomy groupoid}
$Hol(\cF)$ of $\cF$: the objects are points in the manifold and
the arrows are equivalence classes of paths with the same
holonomy.

Since $P_0 A$ is second countable, we can take an open cover $\{
U_i\}$ of $P_0 A$ which consists of countably many small enough
open sets so that in each chart $U_i$, one can choose a
transversal $P_i$ of  the foliation $\cF$. By Proposition 4.8 in
\cite{cf}, each $P_i$ is a smooth manifold with dimension equal to
that of $A$. Let $P=\coprod P_i$, which is a smooth immersed
submanifold of $P_0A$. We can choose $\{U_i\}$ and transversal
$\{P_i\}$ to satisfy the following conditions:
\begin{enumerate}
\item If $U_i$ contains the constant path $0_x$ for some $x\in M$, then
$U_i$ has the transversal $P_i$ containing all constant paths
$0_y$ in $U_i$ for  $y\in M$.
\item If $a(t)\in P_i$ for some $i$, then $a(1-t)\in P_j$ for some
$j$.
\end{enumerate}
It is possible to meet the above two conditions: for (1) we refer
readers to Proposition 4.8 in \cite{cf}. There the result is for
$P_aA$. For $P_0A$, one has to use a smooth reparameterization
$\tau$ with the properties:
\begin{enumerate}
\item $\tau(t)=1$ for all $t \geq 1$ and $\tau(t)=0$ for all $t
  \leq 0$;
\item $\tau'(t)>0$ for all $t\in (0, 1)$.
\end{enumerate}
Then $a^\tau(t):= \tau(t)' a(\tau(t))$ is in $P_0A$ for all $a\in
P_aA$. $\phi_\tau: a\mapsto a^\tau$ defines an injective bounded
linear map from $P_aA \to P_0A$. Therefore, we can adapt the
construction for $P_a A$ to our case by using $\phi_\tau$.
For (2), we define a map $inv:\; P_0 A \to P_0 A $ by
$inv(a(t))=a(1-t)$. Obviously $inv$ is an isomorphism. In
particular, it is open. So we can add $inv(U_i)$ and $inv(P_i)$ to
the collection of open sets and transversals. The new collection
will have the desired property.

 The restriction $Mon(P_0 A)|_P$ of $Mon(P_0 A)$ to $P$ is a finite
dimensional \'etale Lie groupoid\footnote{An \'etale Lie groupoid
is a Lie groupoid such that the source (hence the target) map is a
local diffeomorphism.} \cite{mrw} which we denote by
$\Gamma\underset{\mathbf{t}_1}{\overset{\mathbf{s}_1}{\rightrightarrows}}P$.
For a different transversal $P'$ the restriction of$Mon(P_0 A)$ to
$P'$ is another finite dimensional \'etaleLie groupoid. All these
groupoids are related by ``Morita Equivalence'' which will be
discussed in the next section. One can do the same  to
$Hol(P_0 A)$ and obtain a finite dimensional\'etale Lie groupoid,
which we denote by
$\Gamma^h\underset{\mathbf{t}_1}{\overset{\mathbf{s}_1}{\rightrightarrows}}P$.
Although these groupoids are Morita equivalent to each other, they
are in general not Morita equivalent to the groupoids induced from
$Mon(P_0 A)$. (See also the next section).

We will build a Weinstein groupoid of
$A$ based on this path space $P_0 A$. One can interpret the ``identity section'' as an
embedding obtained from taking constant paths $0_x$, for all $x\in
M$;  the ``inverse'' of a path $a(t)$ as $a(1-t)$;
 the ``source and target maps'' $\bs$ and $\bt$ as taking the end
points of the base path $\gamma(t)$.  According to the two
conditions above, these maps are also well-defined on the finite
dimensional space $P$. Since reparameterizations and projections
are bounded linear operators in Banach space $\ci (I, \R^n)$, the
maps defined above are smooth maps in $P_0 A$, hence in $P$. 

To define the multiplication, notice that for two $A$-paths $a_1$,
$a_0$ in $P_0 A$ such that the base paths satisfy
$\gamma_0(1)=\gamma_1(0)$, one can define a ``concatenation''
\cite{cf}:
\[  a_1\odot a_0 =\left\{
 \begin{array}{rr}
 2a_0(2t), \quad &\mbox{$0\leq t\leq \half$}\\
 2a_1(2t-1),\quad &\mbox{$\half < t\leq 1$}
 \end{array}\right.   \]

Concatenation is a bounded linear operator in the local charts,
hence is a smooth map.  However it is not associative. Moreover it is 
not well-defined on $P$. If we quotient out by the equivalence
relation induced by $\cF$, concatenation is associative and well-defined.
However, after quotient out by the equivalence, we may not end up
with a smooth manifold any more. To overcome the difficulty, our
solution is to pass to the world of differentiable stacks.

\section{Differentiable stacks and Lie groupoids} \label{sect: stack}
 The notion of stacks has been extensively studied in algebraic geometry for
the past few decades (see for
example\cite{dm}\cite{v1}\cite{lmb}\cite{beffgk}). However stacks can also be defined over other categories, such as the
category of topological spaces and category of smooth manifolds
(see for example \cite{SGA4} \cite{pronk} \cite{v2}
\cite{bx}\cite{metzler}). In this section we collect certain facts about stacks in the
differentiable category that will be used in next sections. Many
of those already appeared in literatures (see for example
\cite{pronk} \cite{bx}\cite{metzler}).

\subsection{Definitions}
Let $\cC$ be the category of differentiable manifolds\footnote{For
the future use, we don't require the manifolds to be Hausdorff.}.
A stack over $\cC$ is a category fibred in groupoids satisfying
two conditions: ``isomorphism is a sheaf'', and ``descent datum is
effective''. Morphisms between stacks are just functors between
categories. We refer to \cite{bx} and \cite{metzler} for the
complete definition and only give here an illustrative example.

\begin{ep}\label{mfd}
Given a manifold $M$, one can view it as a stack over $\cC$. Let
$\underline{M}$ be the category where
\[Obj(\underline{M})=\{ (S, u): S\in \cC, u\in Hom(S, M)\}, \]and
a morphism $(S, u) \to (T, v)$ of objects is a morphism $f: S\to
T$ such that $u=v\circ f$. This category\footnote{In a fancier
term, this is the category associated to the functor of points
of $M$.} encodes all information\footnote{cf. Yoneda lemma.} of $M$
in the sense that the morphisms between stacks $\underline{M}$ and
$\underline{M}'$ all come from the ordinary morphisms between $M$
and $M'$. In this way, the notion of stacks generalizes the notion
of manifolds. A stack isomorphic to $\underline{M}$ for some
$M\in\cC$ is called {\bf representable}. In view of this we can
identify $M$ with $\underline{M}$ and treat it as a differentiable
stack.
\end{ep}


We refer the reader to \cite{bx} or \cite{metzler} (Definition 48,
49) for the definition of monomorphisms and epimorphisms of
stacks.

\begin{defi}[representable (surjective) submersions \cite{bx}]
A morphism $f:\cX\to\cY$ of stacks is a  representable
submersion if for every manifold $M$ and every morphism $M\to\cY$
the fibred product $\cX\times_{\cY}M$ is representable and the
induced morphism $\cX\times_{\cY}M\to M$ is a submersion. A
representable submersion is  surjective, if it is furthermore an
epimorphism.
\end{defi}

\begin{defi}[differentiable stacks \cite{bx}] \label{stack}
A  differentiable stack $\cX$ is a stack over $\cC$ together
with a representable surjective submersion $\pi: X\to\cX$ from a
Hausdorff (smooth) manifold $X$. $X$ together with the structure
morphism $\pi: X\to\cX$ is called an atlas for$\cX$.
\end{defi}

\begin{remark}
For a manifold $M$, the stack $\underline{M}$ is by definition
a differentiable stack.
\end{remark}

\begin{ep} \label{bg}
Let $G$ be a Lie group. The set of  principal $G$-bundles forms a 
stack $BG$ in the following way. The objects of $BG$ are
\[ Obj (BG)=\{\pi:  P\to M | P \;\text{ is a principal $G$-bundle over
  $M$}. \} \]
A morphism between two objects $(P, M)$ and $(P', M')$ is a
morphism $M\to M'$ and $G$-equivariant morphism $P\to P'$ covering
$M\to M'$. Moreover $BG$ is a differentiable stack. One can take a point $pt$
to be an atlas. The map $\pi:\underline{pt}\to BG$ defined by
\[(f: M \to pt) \mapsto (M\times_{f,pt, pr} G), \]
(where $pr$ is the projection from $G$ to the point $pt$) is a
representable surjective submersion.
\end{ep}

We have the following two easy properties:

\begin{lemma} [composition]\label{comp}
The composition of two representable (surjective) submersions is
still a representable (surjective) submersion.
\end{lemma}
\begin{proof}
For any manifold $U$ with map $U\to\cZ$, consider the following
diagram
\[
\begin{CD}
\cX\times_\cY \cY \times_\cZ U
@>\tilde{f}>>\cY \times_\cZ U @>\tilde{g}>> U \\
@VVV @VVV @VVV \\
\cX @> f>> \cY @>g>> \cZ .
\end{CD}
\]
Since $f$ and $g$ are representable submersions, $\cY \times_\cZ
U$ is a manifold so that  $\cX\times_\cZ U=\cX\times_\cY\cY
\times_\cZ U$ is  also a manifold. Since $\tilde{g}$ and
$\tilde{f}$ are submersions, $\tilde{g}\circ \tilde{f}$ is also a
submersion. The composition of epimorphisms is still an
epimorphism.
\end{proof}

\begin{lemma}[base change] \label{BC}
In the following diagram
\begin{equation}\label{bc}
\begin{CD}
\cX \times_{\cY}\cZ @>g>> \cZ \\
@VVV @VVV \\
\cX @>f>>\cY,
\end{CD}
\end{equation}
where $\cX$ and $\cY$ are differentiable stacks (but not
necessarily $\cZ$), if $f$ is a representable (surjective)
submersion, then so is $g$.
\end{lemma}
\begin{proof}
For any manifold $U$ mapping to $\cZ$, we have the following
diagram
\[
\begin{CD}
\cX \times_{\cY}\cZ\times_{\cZ} U @>h>> U \\
@VVV @VVV \\
\cX \times_{\cY}\cZ @>g>> \cZ
\end{CD}
\]
Composing the above diagram with \eqref{bc}, one can see that
$\cX{\times}_{\cY}\cZ\times_{\cZ} U=\cX \times_{\cY} U$ is a
manifold and $h$ is a submersion because $f$ is a representable
submersion. Therefore $g$ is a representable submersion. Moreover,
the base change of an epimorphism is clearly still an epimorphism.
\end{proof}
\begin{remark} In general, we call the procedure of obtaining $g$
from $f$ the base change of $\cX \to \cY$ by $\cZ \to \cY$.
\end{remark}

\begin{defi}[smooth morphisms of differentiable stacks]
A morphism $f:\cX\to\cY$ of differentiable stacks is smooth if for
any atlas $g:X\to\cX$ the composition $X\to\cX\to\cY$ satisfies
the following: for any atlas $Y\to\cY$ the induced morphism
$X\times_{\cY}Y\to Y$ is a smooth morphism of manifolds.
\end{defi}

\begin{defi}[immersions, closed immersion, \'etale map, injective immersion. \cite{metzler}]
A morphism $f:\cX\to\cY$ of stacks is an immersion (resp. a closed
immersion, an \'etale map) if for every representable submersion
$M\to\cY$ from a manifold $M$ the product $\cX\times_{\cY}M$ is a
manifold and the induced morphism $\cX\times_{\cY}M\to M$ is an
immersion (resp. a closed immersion, an \'etale map) of manifolds.
An injective immersion is an immersion which is also a
monomorphism.
\end{defi}

A differentiable stack $\cX$ is called {\bf \'etale} if there is
an atlas $\pi: X\to\cX$ with $\pi$ being \'etale.

\begin{lemma}\label{sm}
A morphism $\cX\to\cY$ is smooth if and only if there exist an
atlas $X\to\cX$ of $\cX$ and an atlas $Y\to\cY$ of $\cY$ such that
the induced morphism $X\times_{\cY}Y\to Y$ is a smooth morphism of
manifolds.
\end{lemma}

\begin{proof}
One implication is obvious. Suppose that $X\times_{\cY}Y\to Y$ is
smooth. Let $T\to\cX$ be another atlas. Then using base change of
$X\times_\cY Y\to \cX$,  $T\times_{\cX}X\times_{\cY}Y$ is a
manifold and $T\times_{\cX}X\times_{\cY}Y \to X\times_\cY Y$ is a
submersion, hence a smooth map. The map
$T\times_{\cX}X\times_{\cY}Y\to Y$ factors as
$T\times_{\cX}X\times_{\cY}Y\to X\times_{\cY}Y\to Y$. Hence
$T\times_{\cX}X\times_{\cY}Y\to Y$ is smooth. It also factors as
$T\times_{\cX}X\times_{\cY}Y\to T\times_{\cY}Y\to Y$. Similarly,
the map $T\times_{\cX}X\times_{\cY}Y\to T\times_{\cY}Y$ is a
submersion, hence $T\times_{\cY}Y\to Y$ is smooth.

Now assume that $U\to\cY$ is an atlas of $\cY$. The induced map
$T\times_{\cY}Y\times_{\cY}U\to Y\times_{\cY}U$ is smooth because
it is the base-change of a smooth map $T\times_{\cY}Y\to Y$ by a
submersion $Y\times_{\cY}U\to Y$. One can find a collection of
locally closed submanifolds in $Y\times_{\cY} U$ which form an
open covering family of $U$.
Since being \'etale is a local property, it follows that
$T\times_{\cY}U\to U$ is smooth as well.
\end{proof}

\begin{lemma}\label{immersion}
A morphism from a manifold $X$ to a differentiable stack $\cY$ is
an immersion if and only if $X\times_{\cY} U\to U$ is an immersion
for some atlas $U\to\cY$.
\end{lemma}
\begin{proof}
One implication is obvious. If $X\times_{\cY} U\to U$ is an
immersion, let $T\to\cY$ be any submersion from a manifold $T$.
The map $X\times_{\cY}U\to U$ is transformed by base-change by a
submersion  $U\times_{\cY}T\to U$ to a map
$X\times_{\cY}U\times_{\cY}T\to U\times_{\cY}T$, which is an
immersion since being an immersion is preserved by base-change.
One can find a collection of locally closed submanifolds $\{T_i\}$
in $U\times_{\cY} T$ which forms a family of charts of $T$.
Moreover $X \times_\cY T$ is a manifold because $T$ is an atlas of
$\cY$. Using base changes, one can see that $X\times_{\cY} T_i \to T_i$
is an immersion and that $\{ X\times_\cY T_i\}$ is an open
covering family of $X\times_\cY T$. Since being an immersion is
local property
, it follows that
$X\times_{\cY}T\to T$ is an immersion as well.
\end{proof}

Similarly we have
\begin{lemma}\label{closedimmersion}
A morphism $\cX\to\cY$ of differentiable stacks is a closed
immersion if and only if $\cX\times_{\cY} U\to U$ is a closed
immersion for some atlas $U\to\cY$.
\end{lemma}
\subsection{Stacks v.s. groupoids}
Next we explain the relationship between stacks and groupoids.
\subsubsection{From stacks to groupoids}\label{s-g}
Let $\cX$ be a stack. Pick an atlas $X_0\to\cX$, we can form
\[X_1:=( X\times_{\cX}X)\rightrightarrows X\] with the two maps
being projections from the first and second factors. By the
definition of an atlas, $X_1$ is a manifold and has a natural
groupoid structure with source and target maps the two maps above.
We call this groupoid a presentation of $\cX$. Different atlases
give rise to different presentations (see for example \cite{v1},
Appendix). An \'etale differential stack has a presentation by an
\'etale groupoid.

\begin{ep}
In Example \ref{mfd} we have the stack $\underline{M}$ with the atlas
$M\to\underline{M}$. $M\times_{\underline{M}}M$ is just the
diagonal in $M\times M$, thus is isomorphic to $M$. Hence we have
a groupoid $M\rightrightarrows M$ with two maps both equal to the
identity. This is clearly isomorphic (as a groupoid) to the
transformation groupoid $\{id\}\times M\rightrightarrows M$. Here
$\{id\}$ is the trivial group with one element.
\end{ep}

\begin{ep}
In Example \ref{bg}, the stack $BG$ can be simply presented by a
point $pt$. The fibre product $pt\times_{BG} pt$ equals to $G$. So
the groupoid presenting $BG$ is simply $G\rightrightarrows pt$.
\end{ep}

\subsubsection{From groupoids to stacks}\label{g-s}
Conversely
, given a groupoid
$G_1\underset{\mathbf{t}}{\overset{\mathbf{s}}{\rightrightarrows}}G_0$,
one can associate a (quotient) stack $\cX$ with an atlas $G_0\to\cX$
such that $G_1=G_0\times_{\cX} G_0$. In algebraic geometry, this
result can be found in for example \cite{v1}. Here we recall the
construction given in \cite{bx} for differentiable stacks. We
begin with several well-known definitions.

\begin{defi}[groupoid action]
A Lie groupoid
$G_1\underset{\mathbf{t}}{\overset{\mathbf{s}}{\rightrightarrows}}G_0$
action on a manifold $M$ from the right (or left) consists of the
following data: a moment map $J: M\to G_0$ and a smooth map $\Phi:
M \times_{J,\bt} G_1 \,(\text{or } G_1\times_{\bs, J} M) \to M$
such that
\begin{enumerate}
\item $J(\Phi(m, g))=\bs(g)$  (or $J(\Phi(g, m) =\bt(g)$);
\item $\Phi(\Phi(m, g), h)=\Phi(m, gh)$ (or $\Phi(h, \Phi(g, m))=\Phi(hg, m)$);
\item $\Phi(m, J(m))=m$ (or $\Phi(J(m), m)=m$).

\end{enumerate}
Here we identify $G_0$ as the identity section of $G_1$. From now
on, the action $\Phi$ is denoted by ``$\cdot$'' for simplicity.
\end{defi}

\begin{defi}[groupoid principal bundles, or torsors]
A manifold $P$ is a  (right or left) principal bundle of a
groupoid $H$ over a manifold $S$, if
\begin{enumerate}
\item there is a surjective submersion $\pi: P \to S$;
\item $H_1$  acts (from the right or left) on $P$ freely
and transitively on each fiber of $\pi$;
\item the moment map $J_H: P \to H_0$ is a surjective submersion.
\end{enumerate}
A right principal $H$ bundle is also called $H$-torsor.
\end{defi}
\begin{remark}
Since the action is free and transitive, one can see
that $P/H=S$.  If $H$ action is free,
$P/H$ is a manifold iff $H$ action is proper. Since $S=P/H$ is a
manifold,  we automatically
obtained that the $H$ action is proper for any $H$-principal
bundle.
\end{remark}

Let
$G=$$(G_1\underset{\mathbf{t}}{\overset{\mathbf{s}}{\rightrightarrows}}G_0)$
be a Lie groupoid\footnote{From now on, we also use one letter $G$
(or $H$, etc.) to represent the groupoid
$G_1\underset{\mathbf{t}}{\overset{\mathbf{s}}{\rightrightarrows}}G_0$
(or
$H_1\underset{\mathbf{t}}{\overset{\mathbf{s}}{\rightrightarrows}}H_0$,
etc.).}. We follow the construction in \cite{bx}. Denote by
$BG$ the category of
right $G$-principal bundles. An object $Q$ of $BG$ over $S\in \cC$
is a right $G$-principal bundle over $S$. A morphism between two
$G$ torsors $\pi_1:Q_1\to S_1$ and $\pi_2:Q_2\to S_2$ is a smooth
map $\Psi$ lifting the morphism $\psi$ between the base manifolds
$S_1$ and $S_2$ (i.e. $\psi\circ\pi_1=\Psi\circ\pi_2$) such that
$\Psi$ is $G_1$-equivariant, i.e. $\Psi(q_1 \cdot g) =\Psi(q_1)
\cdot g$ for $(q_1, g) \in M \times_{J_1,\bt} G_1$, where $J_i$
and $\pi_i$ are moment maps and the projections of torsors $Q_i$,
$i=1, 2$.

\begin{note} Above condition implies that $J_2\circ\Psi = J_1$.\end{note}

This makes $BG$ a category over $\cC$. It is a differentiable
stack presented by the Lie groupoid
$G_1\underset{\mathbf{t}}{\overset{\mathbf{s}}{\rightrightarrows}}G_0$:
an atlas $\phi: G_0\to BG$ can be constructed as follows. For $f:
S\to G_0$, we assign a manifold $Q=S \times_{f,\bt} G_1$.  ($Q$ is
a manifold because $\bt$ is a submersion.) The projection $\pi:
Q\to S$ is given by the first projection and the moment map $J: Q
\to G_0$ is the second projection composed with $\bs$. The
groupoid action is defined by
\[ (s, g) \cdot h= (s, gh), \text{ for all} \; (s, g) \in S
\times_{f,\bt} G_1, h \in G_1. \] The $\pi$-fiber is simply a copy
of the $\bt$-fiber, therefore the action of $G_1$ is free and
transitive.

 $\phi$ is a representable surjective submersion
and $G_1= G_0 \times_{\phi,\phi} G_0$ fits in the following
diagram:
\[
\begin{CD}
  G_1 @>\bt>> G_0 \\
  @V \bs VV @V \phi VV \\
  G_0@>\phi>> BG.
\end{CD}
\]
For complete proofs, we refer to literatures more specific on this
subject, for example \cite{pronk}, \cite{bx} and \cite{metzler}.
\begin{ep}
In the case of the trivial transformation groupoid $\{id\}\times
M\rightrightarrows M$ it's easy to see that the stack constructed
above is $\underline{M}$.
\end{ep}
\subsubsection{Morita equivalence}
To further explore the correspondence between stacks and
groupoids, we need the following definition.

\begin{defi}[Morita equivalence \cite{mrw}]
Two Lie groupoids $G$ and $H$ are  Morita equivalent if there
exists a manifold $E$, such that
\begin{enumerate}
\item $G$ and $H$ act on $E$ from the left and right respectively
 with moment maps $J_G$ and $J_H$ and the two actions commute;
\item The moment maps are surjective submersions;
\item The groupoid actions on the fibre of the moment maps are
free and transitive. Such an $E$ is called a {\bf (Morita)
bibundle} of $G$ and $H$.
\end{enumerate}
\end{defi}

\begin{prop}[see \cite{pronk} \cite{bx} \cite{metzler}]
Two Lie groupoids present isomorphic differential stacks
if and only if they are Morita equivalent.
\end{prop}

\subsubsection{1-morphisms}
(1-)morphisms between stacks can be realized on the level of
groupoids.  These had been studied in detail in \cite{hs} and called
Hilsum-Skandalis (HS) morphisms, following \cite{hisk}. Such morphisms are
also called generalized morphisms of Lie groupoids sometimes (see
\cite{moerdijk} and references therein).

\begin{defi}[HS morphisms \cite{hs}]
A   Hilsum-Skandalis (HS) morphism of Lie
groupoids from $G$ to $H$ is a triple $(E, J_G, J_H)$ such that:
\begin{enumerate}
\item The bundle $J_G: E\to G_0$ is a right $H$-principal bundle with moment map $J_H$;
\item $G$ acts on $E$ from the left with moment map $J_G$;
\item The actions of $G$ and $H$ commute, i.e. $(g\cdot x) \cdot
  h=g\cdot (x\cdot h)$.
We call $E$ an {\bf HS bibundle}.
\end{enumerate}
\end{defi}
\begin{remark}
\item{i)} In the above definition, (3) implies that $J_H$ is $G$
  invariant and $J_G$ is $H$ invariant.

\item{ii)} For a homomorphism of Lie groupoids
$f:$$(G_1\underset{\mathbf{t}}{\overset{\mathbf{s}}{\rightrightarrows}}G_0)$$\to$$(H_1\underset{\mathbf{t}}{\overset{\mathbf{s}}{\rightrightarrows}}H_0)$,
one can form an HS morphism via the bibundle $G_0 \times_{f,H_0,
\bt} H_1$ \cite{hm}. Thus the notion of HS morphisms generalizes
the notion of Lie groupoid morphisms.

\item{iii)}  The identity HS morphism of $G_1\underset{\mathbf{t}}{\overset{\mathbf{s}}{\rightrightarrows}}G_0$ is given by
 $G_0 \times_{\bt} G_1 \times_{\bs} G_0$. An HS
  morphism is invertible if the bibundle is not only right
  principal but also left principal. In other words, it is a Morita equivalence.

\item{iv)} Two HS morphisms $E$: $(G_1\rightrightarrows G_0)$$\to$$(H_1\rightrightarrows H_0)$ and
$F$: $(H_1\rightrightarrows H_0)$$\to$$(K_1\rightrightarrows K_0)$
 can be composed to obtain an HS morphism $(G_1\rightrightarrows
G_0)$$\to$$(K_1\rightrightarrows K_0$ ) with bibundle
$E\times_{H_0} F$, where $H_1$ acts by $(x, y)\cdot h=(x h, h^{-1}
y)$ ($G_1$ and $K_1$ still have left-over actions on it).
Composition is not associative (see for example \cite{hm}). 
\end{remark}

\begin{prop}[HS and smooth morphism of stacks]\label{hs}
HS morphisms of Lie groupoids correspond to smooth morphisms of
differentiable stacks. More precisely, an HS morphism $E$:
$(G_1\underset{\mathbf{t}}{\overset{\mathbf{s}}{\rightrightarrows}}G_0)$$\to$$(H_1\underset{\mathbf{t}}{\overset{\mathbf{s}}{\rightrightarrows}}H_0)$
induces a smooth morphism of differentiable stacks $\phi_E: BG \to
BH$. On the other hand, given a smooth morphism $\phi$: $\cX \to
\cY$ and atlases $G_0\to\cX, H_0\to\cY$, $\phi$ induces an HS
morphism $E_\phi$:
$(G_1\underset{\mathbf{t}}{\overset{\mathbf{s}}{\rightrightarrows}}G_0)$$\to$$(H_1\underset{\mathbf{t}}{\overset{\mathbf{s}}{\rightrightarrows}}H_0)$,
where $(G_0\times_{\cX}G_0=G_1)\rightrightarrows G_0$ and
$(H_0\times_{\cY}H_0=H_1)\rightrightarrows H_0$ present $\cX$ and
$\cY$ respectively.
\end{prop}
\begin{proof}
Suppose $(E, J_G, J_H)$ is an HS morphism. Given a right
$G$-principal bundle  $P$ over $S$ with moment map $J_P$, we form
$Q=P\times_{G_0}E/G$, where the $G$-action is given by
\[(p, x)\cdot g=(pg, g^{-1}x), \;\text{if}\;
J_P(p)=J_G(x)=\bt(g).\] Since the action of $G$ is free and proper
on $P$, the $G$-action on $P\times_{G_0}E$ is also free and
proper. So $Q$ is a manifold. In the following steps, we will show
that $Q$ is a $H$-torsor, then we can define $\phi_E$ by
$\phi_E(P)=Q$.

\begin{enumerate}
\item Define $\pi_Q: Q\to S$ by $\pi_Q([(p, x)])=\pi_P (p)$. Since
$\pi_P: P\to S$ is $G$ invariant, $\pi_Q$ is a well-defined
smooth
  map. Since any curve $\gamma(t)$ in $S$ can be pulled back by
  $\pi_P$ as $\tilde{\gamma}(t)$ in $P$, it ($\gamma(t)$) can be pulled back by
  $\pi_Q$ to $Q=P\times_{G_0}E$ as $[(\tilde{\gamma}(t),
    x)]$. Therefore $\pi_Q$ is a surjective submersion.

\item Define $J_Q: Q\to H_0$ by $J_Q([p, x])=J_H(x)$. Since $J_H$
is $G$ invariant, $J_Q$ is well-defined and smooth.

\item Define $H$ action on $Q$ by $[(p, x)]\cdot h=[(p, xh)]$. It
is well defined since $G$ and $H$ actions commute. If $[(p,
x)]\cdot h=[(p, x)]$, then there exist $g\in G_1$, such that $(pg,
g^{-1}xh)=(p, x)$. Since the $G$ action is free on $P$ and the $H$
action is free on $E$, we must have $g=1$ and $h=1$.
Therefore the $H$ action on $Q$ is free.

\item If $[(p, x)]$ and $[(p', x')]$ belong to the same fibre of
$\pi_Q$, i.e. $\pi_P(p)=\pi_P(p')$, then there exists a $g\in
G_1$, such that $p'=pg$. So $[(p, x)]=[(p', g^{-1} x)]$. Since
$J_G(x')=J_P(p')=\bs(g)=J_G(g^{-1}x)$, there exists a $h \in H_1$,
such that $x'h=g^{-1}x$. So $[(p', x')] h=[(p, x)]$, i.e. the
$H$ action on $Q$ is transitive.
\end{enumerate}

On the level of morphisms, we define a map which takes a
morphism of right $G$ principal bundles $f: P_1\to P_2$ to a
morphism of right $H$ principal bundles
\[ \tilde{f}: P_1\times E/G \to P_2 \times E /G, \; \text{given by} \; [(p,
  x)]\mapsto [(f(p), x)]. \]
Therefore $\phi_E$ is a map between stacks. The smoothness of $\phi_E$ follows
from the following claim and Lemma \ref{sm}. \\
{\em Claim}: As manifolds $E\cong H_0 \times_\cY G_0$, where the
map $G_0 \to \cY$ is the composition of the atlas projection
$\pi_G: G_0\to \cX$ and $\phi_E$. The two moment maps $J_H$ and
$J_G$ coincide with the projections from $H_0\times_\cY G_0$ to
$H_0$ and $G_0$ respectively. \\
{\em Proof of the Claim:} Since the inclusion of the category of
manifolds into the category of stacks is faithful, it suffices to
show $E$ and $G_0\times_\cY G_0$ are isomorphic as stacks.

Examining the definition of fibre product of stacks \cite{bx},
we see that an object in $H_0\times_\cY G_0$ over a manifold $S\in
\cC$ is $(f_H, f, f_G)$ where $f_H: S\to
H_0$, $f_G: S\to G_0$ and $f$ is an $H$ equivariant map fitting
inside the following diagram:
\[
\begin{CD}
S\times_{f_H, H_0, \bt} H_1 @>f>> S\times_{f_G, G_0, J_G} E \\
@VVV  @VVV \\
S @>id>> S,
\end{CD}
\]
where we use $(x, e)\mapsto [(x, 1_x, e)]$ to identify the target
of $f$ with $(S\times_{f_G, \bt}G_1\times_{\bs, J_G}E)/G$ which is
the image of the trivial torsor $S\times_{f_G, \bt}G_1$ under the
map $\phi_E\circ \pi_G$. Then by $x\mapsto pr_E\circ f(x, 1_x)$,
$f$ gives a map $\psi_f: S\to E$, which is an object in the stack
$E$.

On the other hand for any $\psi: S\to E$, one can construct a map
$f:S\times_{f_H, H_0, \bt} H_1 \to S\times_{f_G, G_0, J_G} E $ by
$f(x, h)=(x, \psi_f(x)\cdot h)$. Moreover, $f_H$ and $f_G$ are
simply the compositions of $\psi$ with the moment maps of $E$.

It is not hard to verify that this gives an isomorphism between
these two stacks.

Finally, from the construction above, it is not hard to see that
the moment maps are exactly the projections from $H_0\times_\cY
G_0$ to $H_0$ and $G_0$.  \hfill $\btd$

We sketch the proof of the second statement (which is not used in
the remaining content of this paper). We have morphisms
$G_0\to\cX\overset{\phi}{\lra}\cY$ and $H_0\to\cY$. Take the
bibundle $E_\phi$ to be $G_0\times_{\cY} H_0$. It's not hard to
check that $E_\phi$ satisfies required properties.
\end{proof}
\begin{remark} A different form of this proposition can be found
in \cite{pronk}, in which a proof from bicategorical viewpoint was
given.
\end{remark}

In view of Proposition \ref {hs}, the fact that the composition of
HS morphisms is not associative can be understood by the fact that
compositions of 1-morphisms of stacks are associative up to
2-morphisms of stacks.

\subsubsection{2-morphisms}
As morphisms in differentiable stacks corresponds to HS morphisms,
2-morphisms also have their exact correspondence in the language
of Lie groupoids. Recall that morphisms of stacks are
functors between categories, and a 2-morphism of stacks between
two morphisms is a natural transformation between these two
morphisms viewed as functors. We have 2-morphisms of groupoids
defined as following:
\begin{defi}[2-morphisms \cite{pronk}\cite{metzler}]
Let $(E^i, J^i_G, J^i_H)$ be two HS morphisms from the Lie
groupoid $G$ to $H$. A  2-morphism from $(E^1, J^1_G, J^1_H)$
to $(E^2, J^2_G, J^2_H)$ is a bi-invariant isomorphism from $E^1$
to $E^2$.
\end{defi}
\begin{remark}\label{2morp}
\item{i)}  If the two HS morphisms are given by groupoid
homomorphisms $f$ and $g$ between $G$ and $H$, then a 2-morphism
from $f$ to $g$ is just a smooth map $\alpha: G_0 \to H_1$ so that
$f(x)=g(x)\cdot \alpha(x) $ and $\alpha(\gamma x)=g(\gamma) \alpha
(x) f(\gamma)^{-1}$, where $x\in G_0$ and $\gamma\in G_1$. So it
is easy to see that not every two morphisms can be connected by a
2-morphism and when they do, the 2-morphism may not be unique (for
example this happens when the isotropy group is nontrivial and
abelian).

\item{ii)} From the proof of Proposition \ref{hs}, one can see
that a 2-morphism between HS morphisms corresponds to a 2-morphism
between the corresponding (1-)morphisms on the level of stacks.
\end{remark}

\subsection{Fibre products and submersions}

Invariant maps are a convenient way to produce maps between stacks
that we will use later in the construction of the WEinstein groupoids.

\begin{lemma} \label{im} 
Given a Lie groupoid
$G_1\underset{\mathbf{t}}{\overset{\mathbf{s}}{\rightrightarrows}}G_0$
and a manifold $M$, any $G$-invariant map $f: G_0\to M$ induces a
morphism between differentiable stacks $\bar{f}: BG \to M$ such
that $f=\bar{f}\circ \phi$, where $\phi: G_0 \to BG$ is the
covering map of atlases.
\end{lemma}
\begin{proof}
Since $f$ is $G$ invariant, $f$ introduces a morphism between Lie
groupoids:
$(G_1\underset{\mathbf{t}}{\overset{\mathbf{s}}{\rightrightarrows}}G_0)$$\to
(M\rightrightarrows M)$. By Proposition \ref{hs} it gives a smooth
morphism between differentiable stacks. More precisely, let $Q\to
S$ be a $G_1$ torsor over $S$ with moment map $J_1$ and projection
$\pi_1$. Since the  $G$ action on the  $\pi_1$-fibre is free and
transitive, we have $S=Q\times_{f\circ J_1,id} M /G_1$. Notice
that a $(M\rightrightarrows M$)-torsor is simply a manifold $S$
with a smooth map to $M$. Then $\bar{f}(Q)$ is the morphism $J_2:
S\to M$ given by $J_2(s)=f\circ J_1(q)$, where $q$ is any preimage
of $s$ by $\pi$ (it is well defined since $f$ is $G$-invariant).
For any map $a: S\to G_0$, the image under $\phi$ is $Q_a=S
\times_{a,\bt} G_1$, and $\bar{f}(Q_a)$ is the map $f\circ a$
since $f$ is $G$-invariant. Therefore $f=\bar{f}\circ \phi$.
\end{proof}

\begin{defi}[(surjective) submersions]
A morphism $f:\cX\to\cY$ of differentiable stacks is called a 
submersion\footnote{This is different from the definition in
\cite{metzler}.} if for any atlas $M\to\cX$, the composition
$M\to\cX\to\cY$ satisfies the following: for any atlas $N\to\cY$
the induced morphism $M\times_{\cY}N\to N$ is a submersion. A 
surjective submersion is a submersion which is also an
epimorphism.
\end{defi}

\begin{remark}
In particular, a representable submersion is a submersion. But the
converse is not true: for example the source and target maps $bbs$ and
$bbt$ that we will define in the next section are submersions but not
representable submersions in general.
\end{remark}

Here we introduce submersions instead of representable submersions
mainly because we will use later the following result about fibred
products:

\begin{prop}[fibred products]\label{fp}
Let $Z$ be a manifold and $f:\cX\to Z$ and $g:\cY\to Z$ be two morphisms of differentiable
stacks. If either $f$ or $g$ is a submersion, then
$\cX \times_Z\cY$ is a differentiable stack.
\end{prop}
\begin{proof}
Assume that $f:\cX \to Z$ is a submersion. By definition, for any
atlas  $X\to\cX$, the composition $X\to\cX\to Z$ is a submersion.
Let $Y$ be a presentation of $\cY$, then $X\times_Z Y$ is a
manifold. To see that $\cX \times_Z \cY$ is a differentiable
stack, it suffices to show that there exists a representable
surjective submersion from $X\times_Z Y$ to $\cX\times_Z \cY$. By
Lemma \ref{BC}, $X\times_Z Y \to \cX\times_Z Y$ and $\cX\times_Z Y
\to\cX\times_Z\cY$ are representable surjective submersions. By
Lemma \ref{comp}, their composition is also a representable
surjective submersion.
\end{proof}

\begin{lemma}\label{pa}
Let $\cX, \cY$ be stacks with maps $\cX\to Z$ and $\cY\to Z$ to a
manifold $Z$, one of which a submersion, and let $X\to\cX,
Y\to\cY$ be atlases for $\cX$ and $\cY$ respectively. Then
$X\times_Z Y\to\cX\times_Z\cY$ is an atlas for $\cX\times_Z\cY$.
\end{lemma}

\begin{proof}
Note that $X\times_Z Y$ is a manifold because one of $\cX\to Z$
and $\cY\to Z$ is a submersion. $X\times_Z Y\to\cX\times_Z\cY$
factors into $X\times_Z Y\to\cX\times_Z Y\to\cX\times_Z\cY$.
$X\times_Z Y\to\cX\times_Z Y$ is a representable surjective
submersion because $X\to\cX$ is. $\cX\times_Z Y\to\cX\times_Z\cY$
is a representable surjective submersion because $Y\to\cY$ is.
Thus $X\times_Z Y\to\cX\times_Z\cY$ is a representable surjective
submersion.
\end{proof}

\begin{ep}\label{pa2}
In the situation of Lemma \ref{pa}, put $X_1=X\times_{\cX}X$ and
$Y_1=Y\times_{\cY} Y$, then $\cX\times_Z\cY$ is presented by the
groupoid $$(X_1\times_Z Y_1\rightrightarrows X\times_Z Y).$$ This
follows from the fact that $$(X\times_Z
Y)\times_{\cX\times_Z\cY}(X\times_Z Y)\cong
(X\times_{\cX}X)\times_Z (Y\times_{\cY} Y).$$
\end{ep}

\begin{lemma}\label{ssts} 
If a $G$ invariant map $f: G_0\to M$ is a submersion, then the
induced map $\bar{f}: BG \to M $ is a submersion of differentiable
stacks.
\end{lemma}

\begin{proof}
Let $U\to M$ be a morphism of manifolds. Using the base change of
the representable surjective submersion $G_0 \to BG$ by the
projection $BG\times_{\bar{f}, M} U \to BG$, we can see that
$BG\times_M U$ is a differentiable stack with atlas $G_0 \times_M
U$. Note that the composition $G_0\times_M U\to BG\times_M U\to U$
is a submersion because it is the base change of $f: G_0\to M$ by
$U\to M$. Now take an atlas $V\to BG\times_M U$ which is a
representable surjective submersion.
$$\begin{diagram}
\node{G_0\times_M U\times_{BG\times_M U}V}\arrow[2]{e}\arrow{se}\node[2]{G_0\times_M U}\arrow[2]{s}\arrow{se}\arrow{ese}\\
\node[2]{V}\arrow[2]{e}\node[2]{BG\times_M U}\arrow{e}\arrow[2]{s}\node{U}\arrow[2]{s}\\
\node[3]{G_0}\arrow{se}\arrow{ese}\\
\node[4]{BG}\arrow{e}\node{M}.
\end{diagram}
$$
We see that $G_0\times_M U\times_{BG\times_M U} V$ is a manifold
and the projections to $G_0\times_M U$ and $V$ are submersions.
The composition $$G_0\times_M U\times_{BG\times_M U} V\to V\to U$$
coincides with
$$G_0\times_M U\times_{BG\times_M U} V\to G_0\times_M U\to
BG\times_M U\to U, $$ which is a surjective submersion.
Hence $V\to U$ is
a submersion.
\end{proof}
 It's not hard to see that the construction of stacks in the category of
smooth manifolds can be extended to the category of Banach
manifolds, yielding the notion of Banach stacks. Many properties
of differentiable stacks, including those discussed here, are
shared by Banach stacks as well. Also, the category of
differentiable stacks can be obtained from the category of Banach
stacks by restricting the base category.

\section{The Weinstein groupoids of Lie algebroids} \label{w}
\subsection{The construction}
Recall that in Section \ref{apath},  given a Lie algebroid $A$, we constructed an \'etale groupoid
$\Gamma\underset{\mathbf{t}_1}{\overset{\mathbf{s}_1}{\rightrightarrows}}P$.
We obtain an \'etale differential stack$\cG(A)$ presented by
$\Gamma\underset{\mathbf{t}_1}{\overset{\mathbf{s}_1}{\rightrightarrows}}P$.
For a different transversal $P'$, the restriction $\Gamma'=Mon(P_0
A)|_{P'}$ is Morita equivalent to $\Gamma$ through the finite
dimensional bibundle $\bs_M^{-1}(P) \cap \bt_M^{-1}(P')$. As we
have seen, this implies that they represent isomorphic
differential stacks. Therefore, we might base our discussion on
$\Gamma\underset{\mathbf{t}_1}{\overset{\mathbf{s}_1}{\rightrightarrows}}P$.

Since $Mon(P_0 A)\rightrightarrows P_0 A $ is Morita equivalent to
$\Gamma\underset{\mathbf{t}_1}{\overset{\mathbf{s}_1}{\rightrightarrows}}P$
through the Banach bibundle $\bs_M^{-1}(P)$, $\cG(A)$ can also be
presented by $Mon(P_0A)$ as a Banach stack.


In this section, we will construct two Weinstein groupoids
$\cG(A)$ and $\cH(A)$ for every Lie algebroid $A$ and  prove
Theorem \ref{integ}.

We begin with $\cG(A)$. We first define the inverse, identity
section, source and target maps on the level of groupoids.

\begin{defi}
Define
\begin{itemize}
  \item $i:$ $(\Gamma\underset{\mathbf{t}_1}{\overset{\mathbf{s}_1}{\rightrightarrows}}P)$$\to$$(\Gamma\underset{\mathbf{t}_1}{\overset{\mathbf{s}_1}{\rightrightarrows}}P)$ by
  $g=[a(\epsilon, t)]\mapsto [a(\epsilon, 1-t)]$, where $[\cdot]$ denotes the homotopy class in $Mon(P_0A)$;
  \item $e: M \to$$(\Gamma\underset{\mathbf{t}_1}{\overset{\mathbf{s}_1}{\rightrightarrows}}P)$ by $x\mapsto 1_{0_x}$, where
  $1_{0_x}$ denotes the identity homotopy of the constant path
  $0_x$;
  \item $\bs:$$(\Gamma\underset{\mathbf{t}_1}{\overset{\mathbf{s}_1}{\rightrightarrows}}P)$$\to M$ by $g=[a(\epsilon,
  t)]\mapsto \gamma(0, 0)(=\gamma(\epsilon, 0),\, \forall \epsilon)$,
  where $\gamma$ is the base path of $a$;
  \item $\bt:$$(\Gamma\underset{\mathbf{t}_1}{\overset{\mathbf{s}_1}{\rightrightarrows}}P)$$\to M$ by $g=[a(\epsilon,
  t)]\mapsto \gamma(0, 1)(=\gamma(\epsilon, 1),\, \forall \epsilon)$;
\end{itemize}
These maps can be defined similarly on
$Mon(P_0A)\rightrightarrows P_0A $. These maps are all bounded
linear maps in the local charts of $Mon(P_0A)$. Therefore they are
smooth homomorphisms between Lie groupoids. Hence, they defined
smooth morphisms between differentiable stacks. We denote the maps
corresponding to $i$, $e$, $\bs$, $\bt$ on the stack level by
$\bar{i}$, $\bar{e}$, $\bar{\bs}$ and $\bar{\bt}$ respectively.
\end{defi}

\begin{lemma}\label{e-immersion}
The maps $\bbs$ and $\bbt$ are surjective submersions. The map
$\bar{e}: M\to \cG(A)$ is an injective immersion. The map
$\bar{i}$ is an isomorphism.
\end{lemma}
\begin{proof}
$\bs$ and $\bt$ restricted to $P$ are $\Gamma$-invariant
and  submersions because any path through $x$ in $M$ can
be lifted to a path in $P$ passing through any given preimage of
$x$. According to Lemma \ref{im} and \ref{ssts}, the induced maps
$\bar{\bs}$ and $\bbt$ are  submersions.

Denote by $e_0$ the restricted map of $e$ on the level of objects:
$e_0: M \to P$. Notice that $e_0$ fits into the following diagram
(which is not  commutative):
\begin{equation}\label{eq: immersion}\begin{diagram}
\node{M\times_{\cG(A)}P}\arrow{e,t}{pr_2}\arrow{s,l}{pr_1}
\node{P}\arrow{s,r}{\pi} \\
\node{M}\arrow{e,t}{\bar{e}}\arrow{ne,t}{e_0}
\node{\cG(A)}\end{diagram}
\end{equation}
Consider $x=(f:U\to M)\in M$, $\bar{e}(x)=U\times_{e_0\circ f,
G_0} G_1$ as a $G$-torsor, and $e_0(x)=(e_0\circ f:U\to G_0)\in
G_0$. Consider also $y=(g: U\to G_0)\in G_0$,
$\pi(y)=U\times_{g, G_0} G_1$. A typical object of $M_i
\times_{\cG} G_0$ is $(x, \eta, y)$ where $\eta$ is a morphism for
$G$-torsors from $\bar{e}(x)$ to $\pi(y)$ over $id_U$ of $U$. Then
by the equivariancy of $\eta$, we have a map $\phi$: $U\to G_1$,
such that $e_0\circ f=g\cdot \phi$. Therefore, we have a map
$\alpha: M\times_{\cG(A)} G_0\to G_1$ given by $\alpha (x, \eta,
y)=\phi$, such that \[e_0 \circ pr_1 =pr_2 \cdot \alpha.\] Since
$\pi$ is \'etale, so is $pr_1$. Moreover, since $e_0$ is an
embedding, $pr_2$ must be an immersion. Therefore, by Lemma
\ref{immersion}, $\bar{e}$ is an immersion.

As $\bs \circ e = \bt \circ e = id$ on the level of groupoids, the
same identity passes to identity on the level of differentiable
stacks. Since $\bbs \circ \bar{e}=\bbt \circ \bar{e}=id$, it is
easy to see that $\bar{e}$ must be monomorphic and $\bbs$ (and
$\bbt$) must be epimorphic.

The map $i$ is an isomorphism of groupoids, hence it induces an
isomorphism at the level of stacks.
\end{proof}

Now we define the multiplication in the infinite dimensional
presentation. First we extend ``concatenation'' to $Mon(P_0A)$.
Consider two elements $g_1, g_0 \in Mon(P_0 A)$ whose base paths
on $M$ are connected at the end points. Suppose $g_i$ is
represented by $a_i(\epsilon, t)$. Define
\[ g_1\odot g_0=[a_1(\epsilon, t) \odot_t a_0
  (\epsilon , t)], \]
where $\odot_t$ means concatenation with respect to  the parameter
$t$ and the $[\cdot]$ denotes the equivalence class of homotopies.

Notice that $\bs\circ\bs_M =\bs\circ\bt_M$ and $\bt\circ\bs_M =
\bt\circ\bt_M$ are surjective submersions by reasoning similar to
that in above. Hence by Lemma \ref{pa} and Example \ref{pa2},
\[Mon(P_0A)\times_{\bs\circ\bs_M, M, \bt\circ\bt_M}Mon(P_0A)\rightrightarrows
P_0A \] with source and target maps $\bs_M \times \bs_M$ and
$\bt_M\times\bt_M$ is a Lie groupoid and it presents the stack
$\cG
\times_{\bbs,M, \bbt} \cG$.\\

Finally let $m$ to be the following smooth homomorphism between
Lie groupoids:
$$\label{multbanach}
\xymatrix@=50pt{Mon(P_0 A)\mathop\times\limits_{\bs\circ\bs_M, M, \bt\circ\bt_M}Mon(P_0 A)\ar[d]^{\bs_M\times \bs_M}\ar@<-1ex>[d]_{\bt_M\times\bt_M}\ar[r]^-{\odot}& Mon(P_0 A)\ar[d]^{\bs_M}\ar@<-1ex>[d]_{\bt_M} \\
P_0 A \times P_0 A\ar[r]^{\odot}& P_0 A \\
}
$$

Multiplication is less obvious for the \'etale presentation
$\Gamma\rightrightarrows P$.  We will have to define the multiplication through an HS
morphism.

Viewing $P$ as a submanifold of $P_0A$, let $E=\bs_M^{-1}(P)\cap
\bt_M^{-1}(m(P\times_M P))\subset Mon(P_0A)$. Since $\bs_M$ and
$\bt_M$ are surjective submersions and $m(P\times_M P)\cong P\times_M
P$ is a submanifold of $P_0A$, $E$ is a smooth manifold.
Since $P$ is a transversal, $\bt_M: E\to m(P\times_M P)$ is
\'etale. Moreover $\dim m(P\times_M P)= 2\dim P -\dim M$.  So
$E$ is finite dimensional. Further notice that $m: P_0A\times
P_0A \to P_0A $ is injective and its ``inverse'' $m^{-1}$ defined on
the image of $m$ is given by
\[ m^{-1}: b(t)\mapsto (b(2t_1), b(1-2t_2)) \;\;\; t_1\in [0,
\frac{1}{2}],\; t_2\in [\frac{1}{2}, 1] \] which is bounded
linear in a local chart. Let $\pi_1=m^{-1}\circ \bt_M: E\to P\times_M P$ and
$\pi_2=\bs_M: E\to P $. Then it is routine to check that $(E,
\pi_1, \pi_2)$ is an HS morphism from $\Gamma\times_M \Gamma
\rightrightarrows P\times_M P$ to $\Gamma \rightrightarrows P$. It
is not hard to verify that on the level of stacks $(E, \pi_1, \pi_2)$ and $m$ give two
1-morphisms  differed by a 2-morphism. Thus,
after modifying $E$ by this 2-morphism, we get another HS-morphism
$(E_m, \pi'_1, pi'_2)$ which presents the same map as $m$.
Moreover, $E_m\cong E$ as bibundles.

Therefore, we have the following definition:

\begin{defi} Define $\bm: \cG(A)\times_{\bbs,\bbt} \cG(A) \to \cG(A)$  to be
the smooth morphism between \'etale stacks presented by $(E_m,
\pi'_1, \pi'_2) $.
\end{defi}

\begin{remark}
If we use $Mon(P_0A)$ as the presentation, $\bm$ is also presented
by $m$.
\end{remark}

\begin{lemma} The multiplication $\bm: \cG(A) \times \cG(A) \to \cG(A)$ is
a smooth morphism between \'etale stacks and is associative up to a 2-morphism,
that is, diagram
$$\begin{diagram}
\node{\cG(A)\mathop\times\limits_{\bs,\bt}\cG(A)\mathop\times\limits_{\bs,\bt}\cG(A)}\arrow{e,t}{id\times\bar{m}}\arrow{s,b}{\bar{m}\times id}\node{\cG(A)\mathop\times\limits_{\bs,\bt}\cG(A)}\arrow{s,r}{\bar{m}}\\
\node{\cG(A)\mathop\times\limits_{\bs,\bt}\cG(A)}\arrow{e,t}{\bar{m}}\node{\cG(A)}
\end{diagram}
$$
is 2-commutative, i.e. there exists a 2-morphism $\alpha: \bm\circ(\bm\times id)
\rightarrow \bm\circ(id\times \bm)$.
\end{lemma}
\begin{proof}
We will establish  the 2-morphism on the level of Banach stacks.
Notice that a smooth morphism in the category of Banach manifolds
between finite dimensional manifolds is a smooth morphism in the
category of finite dimensional smooth manifolds. Therefore, the
2-morphism we will establish gives a 2-morphism for the \'etale stacks.

Take the Banach presentation
$Mon(P_0 A)$, then $\bm$ can simply be presented as a
homomorphism between groupoids as in \eqref{multbanach}. According to
Remark \ref{2morp}, we now construct a 2-morphism $\alpha:
P_0 A\times_M P_0 A \times_M P_0 A \to Mon(P_0 A)$ in the following
diagram
\[
\xymatrix@=50pt{Mon(P_0 A)\mathop\times\limits_{M}Mon(P_0 A)\mathop\times\limits_{M}Mon(P_0 A)\ar[r]^-{m \circ (m \times id)}_-{m \circ (id \times m)}\ar[d]^{\bs_M\times\bs_M\times\bs_M}\ar@<-1ex>[d]_{\bt_M\times\bt_M\times\bt_M}& Mon(P_0 A)\ar[d]^{\bs_M}\ar@<-1ex>[d]_{\bt_M}\\
P_0 A\times_M P_0 A \times_M P_0 A \ar[r]& P_0 A
}
\]

Let $\alpha(a_1, a_2, a_3)$ be the natural rescaling between $a_1
\odot (a_2\odot a_3)$ and $(a_1\odot a_2)\odot a_3$.
Namely, $\alpha(a_1, a_2, a_3)$ is the homotopy class
represented by \begin{equation}\label{rp}
a(\epsilon, t) =
((1-\epsilon)+\epsilon\sigma'(t))a((1-\epsilon)t+\epsilon\sigma(t)),
\end{equation}
where $\sigma(t)$ is a smooth reparameterization such that
$\sigma(1/4)=1/2, \; \sigma(1/2)=3/4.$
In local charts, $\alpha$ is a bounded linear operator. Therefore
it is a smooth morphism between Banach spaces. Moreover, $ m \circ
(m \times id)=  m \circ (id \times m)\cdot\alpha $. Therefore
$\alpha$ serves as the desired 2-morphism.
\end{proof}

One might be curious about whether there are further obstructions to
associativity. There are six ways to multiply four elements in
$\cG(A)$. Put
\[
\begin{split}
F_1&=\bm \circ \bm \times id \circ \bm \times id \times id, \\
F_2&=\bm \circ id \times \bm \circ \bm \times id \times id, \\
F_3&=\bm \circ \bm \times id \circ id \times id \times \bm, \\
F_4&=\bm \circ id \times \bm \circ id \times id \times \bm, \\
F_5&=\bm \circ id \times \bm \circ id \times \bm \times id, \\
F_6&=\bm \circ \bm \times id \circ id \times \bm \times id.
\end{split}
\]
These morphisms fit into the following commutative cube.
$$
\xymatrix@=5pt{
     & & \cG(A)\mathop\times\limits_{M}\cG(A)\mathop\times\limits_{M}\cG(A) \ar[dr]^{\bar{m}\times id} \ar[ddd]^{id\times\bar{m}}& \\
  \cG(A)\mathop\times\limits_{M}\cG(A)\mathop\times\limits_{M}\cG(A)\mathop\times\limits_{M}\cG(A) \ar[urr]^{id\times id\times\bar{m}} \ar[dr]^{\bar{m}\times id\times id} \ar[ddd]_{id\times\bar{m}\times id} & & & \cG(A)\mathop\times\limits_{M}\cG(A)\ar[ddd]^{\bar{m}} \\
     & \cG(A)\mathop\times\limits_{M}\cG(A)\mathop\times\limits_{M}\cG(A)\ar[urr]^{id\times\bar{m}} \ar[ddd]^{\bar{m}\times id} & & & \\
     & & \cG(A)\mathop\times\limits_{M}\cG(A) \ar[dr]^{\bar{m}} & \\
  \cG(A)\mathop\times\limits_{M}\cG(A)\mathop\times\limits_{M}\cG(A) \ar[urr]^{id\times\bar{m}} \ar[dr]^{\bar{m}\times id} & & & \cG(A) \\
     & \cG(A)\mathop\times\limits_{M}\cG(A) \ar[urr]^{\bar{m}} & &
       }
$$
There is a 2-morphism on each face of the cube to connect $F_i$
and $F_{i+1}$ ($F_7=F_1$), constructed as in the last lemma. Let
$\alpha_i: F_i\to F_{i+1}$. Will the composition $\alpha_6 \circ
\alpha_6 \circ ... \circ \alpha_1$ be the identity 2-morphism? If
so, given any two different ways of multiplying four (hence any
number of) elements, different methods to obtain 2-morphisms
between them will give rise to the same 2-morphism. Since
2-morphisms between two 1-morphisms are not unique if our
differential stacks are not honest manifolds,  it is necessary to
study the existence of further obstructions.

\begin{prop} \label{3asso} There is no further obstruction for
associativity of $\bm$ in $\cG(A)$.
\end{prop}
\begin{proof} In the presentation $Mon(P_0A)$ of $\cG(A)$,
the $\alpha_i$'s constructed above can be explicitly expressed as
a smooth morphism: $P_0A\times_MP_0A\times_M P_0A\times_M P_0A \to
Mon(P_0A)$. More precisely,according to the Lemma above,
$\alpha_i(a_1, a_2, a_3, a_4)$ is the natural rescaling between
$F_i(a_1, a_2, a_3, a_4)$ and$F_{i+1}(a_1, a_2, a_3, a_4)$. Here
by abuse of notations, we denote the homomorphism on the groupoid
level also by $F_i$. It is not hard to see that $\alpha_6 \circ
\alpha_6 \circ ... \circ\alpha_1$ is represented by a rescaling
that is homotopic to the identity homotopy between $A_0$-paths.

Therefore, the composed 2-morphism is actually identity since
$Mon(P_0A)$ is made up by the homotopy of homotopy of $A_0$-paths.
We also notice that identity morphism in the category of Banach
manifolds between two finite dimensional manifolds is identity morphism in the
category of finite dimensional smooth manifolds. Therefore, there
is no further obstructions even for 2-morphisms of \'etale stacks.
\end{proof}

Now to show $\cG(A)$ is a Weinstein groupoid, it remains to show
that the identities in item (4) and (5) in Definition \ref{wgpd} hold and
the 2-morphisms in these identities are identity 2-morphisms when restricted to $M$.
Notice that for any $A_0$-path $a(t)$, we have
\[ a(t) \odot_t 1_{\gamma(0)} \sim a(t), \;\; a(1-t) \odot_t a(t)
\sim \gamma(0) ,\] where $\gamma$ is the base path of $a(t)$.
Using i) in Remark \ref{2morp}, we can see that on the groupoid
level $m\circ ((e\circ \bt) \times id)$ and $id$ only differ by a
2-morphism, and the same for the pairs  $m\circ (i \times id)$and
$e\circ \bs$, $\bs \circ m$ and $\bs \circ pr_1$. Therefore the
corresponding identities hold on the level of differentiable
stacks. Transform them to stacks, the rest of the identities also
follow. Moreover, the 2-morphisms (in all presentations of
$\cG(A)$ we have described above) are formed by rescalings. When
they restrict to constant paths in $M$, they are just $id$.

Summing up what we have discussed above, $\cG(A)$ with all the
structures we have given is a Weinstein groupoid over $M$.

We further comment that one can construct another natural
Weinstein groupoid $\cH(A)$ associated to $A$ exactly in the same
way as $\cG(A)$ by the Lie groupoid $Hol(P_0 A)$ or
$\Gamma^h\underset{\mathbf{t}_1}{\overset{\mathbf{s}_1}{\rightrightarrows}}P$
since they are Morita equivalent by a similar reason as their
monodromy counterparts. One can establish the identity section,
the inverse, etc., even the multiplication in exactly the same
way. One only has to notice that in the construction of the
multiplication, the 2-morphism in the associativity diagram is the
holonomy class (instead of homotopy class) of the
reparameterization \eqref{rp}. One can do so because homotopic
paths have the same holonomy. Moreover, by the same reason, there
is no further obstructions for the multiplication on $\cH(A)$.

Finally, we want to comment about the Hausdorffness of the source
fibres (hence the target fibres by the inverse) of $\cG(A)$ and $\cH(A)$.

\begin{defi} An \'etale differentiable stack $\cX$ is Hausdorff
iff the diagonal map\[\Delta: \cX \to \cX \times \cX ,\] is an
closed immersion.
\end{defi}
\begin{remark} In the case when $\cX$ is a manifold, the diagonal
map being a closed immersion is equivalent to its image being
closed. Hence this notion coincides with the usual Hausdorffness
for manifolds.
\end{remark}

Unlike the case of Lie groupoids, the source fibre of $\cG(A)$
or $\cH(A)$ is in general not Hausdorff. (see Example
\ref{nonhaus}). The obstruction lies inside the foliation $\cF$
defined in Section \ref{apath}.

\begin{prop} The source fibre of $\cG(A)$ and $\cH(A)$ is
Hausdorff iff the leaves of the foliation $\cF$ are closed.
\end{prop}
\begin{proof}
We prove this for $\cG(A)$. The proof for $\cH(A)$ is similar. Let $P$
be the \'etale atlas we have chosen. Then the source fibre
$\bbs^{-1}(x)=x\times_{M, \bbs} \cG(A)$ is a differentiable stack
presented by $\bs^{-1}(x)$ by Proposition \ref{fp}. Consider the
following diagram,
\[
\begin{CD}
\bbs^{-1}(x)\times_{\bbs^{-1}(x) \times \bbs^{-1}(x)} \bs^{-1}(x)\times \bs^{-1}(x) @>\delta>> \bs^{-1}(x) \times \bs^{-1}(x) \\
@VVV                  @VVV \\
\bbs^{-1}(x)@>\Delta>> \bbs^{-1}(x) \times \bbs^{-1}(x).
\end{CD}
\]
By a similar argument as in Section \ref{g-s},
$\bbs^{-1}(x)\times_{\bbs^{-1}(x) \times \bbs^{-1}(x)}
\bs^{-1}(x)\times \bs^{-1}(x)$ is isomorphic to
$\Gamma|_{\bs^{-1}(x)}$ and $\delta$ is just $\bs_1\times \bt_1$.
Obviously $\bs_1\times\bt_1$ is an immersion since $\Gamma$ is an
\'etale groupoid. Moreover, the image of $\delta$ is closed by the
following argument: take a convergent sequence $(a^i_0(t),
a^i_1(t))$ of $A_0$-path with the limit $(a_0(t), a_1(t))$.
Suppose that $(a^i_0(t), a^i_1(t))$ is inside the image of
$\delta$, i.e. $a^i_0(t)\sim a^i_1(t))$. Let $\bar{a}$ denote the
inverse path of $a$, we have $\bar{a}^i_0(t) \odot a^i_1(t) \sim
1_x$, i.e. they stay in the same leaf of the foliation $\cF$.
Hence the limit path $\bar{a}_0(t) \odot a_1(t)\sim 1_x$ (i.e.
$(a_0(t), a_1(t))$ is also inside the image of $\delta$) iff the
leaves of $\cF$ are closed.
\end{proof}

\begin{ep}[Non-Hausdorff source-fibres] \label{nonhaus}
Let $M$ be $S^2\times S^2$ with 2-form $\Omega=(\omega,
\sqrt{2}\omega)$. Let the Lie algebroid $A$ over $M$ be $TM\times \R$ with Lie bracket
\[ [ (V, f), (W, g)] = ([V, W], L_V(g) -L_W(h) + \Omega (V, W)),  \]
and anchor the projection onto $TM$ (see \cite{cw} Chapter 16 or \cite{am1}). Let $(a(\epsilon, t), u(\epsilon,
t))$ be an  $A_0$-homotopy, where the first component is in $TM$ and the
second component is in the trivial bundle $\R$. The condition of being
an $A_0$-path here is equivalent to $a=\frac{d}{dt}\gamma$ and
boundary conditions, where
$\gamma$ is the base path.  Moreover, the first component of
equation \eqref{eq-homotopy} is the usual $A_0$-homotopy equation  for $TM$, which
simply induces the homotopy of the base paths. The second component of
equation \eqref{eq-homotopy} is \[\partial_t v - \partial_{\epsilon} u =
\Omega(a, b), \] where $b$ in equation \eqref{eq-homotopy} is $(b, v)$
above. Hence $b=\frac{d}{d\epsilon}\gamma$. Integrate the above
equation and use the boundary condition of $v$, we have
\[\int_0^1 u(0, t) dt- \int_0^1 u(1, t) dt= \int_{\gamma} \Omega. \]
 Let the period group $\Lambda$ of $\Omega$ at a point
$x\in M$ be
\[ \Lambda_x =\int_\gamma \Omega, \quad [\gamma]\in \pi_2(M, x). \]
Since $M$ is simply connected, one can actually show that $(\gamma(0, t), u(0, t)) \sim (\gamma(1,
t), u(1, t))$ iff $\gamma(0, t)$ and $\gamma(1, t)$ have the same end
points and $\int_0^1 (u_0-u_1) dt \in \Lambda$.
Then in this case, since $\sqrt{2}$ is irrational, $\Lambda_x$ is dense
in $\R$ for all $x$. Hence, there exist sequences $u_0^i\to u_0$ and
$u_1^i\to u_1$
such that $\int_0^1 (u_0^i-u_1^i) \in \Lambda$ but the limit $\int_0^1
(u_0-u_1) \notin \Lambda$. Hence the leaves of the foliation $\cF$ are
not closed in this case.
\end{ep}
\subsection{The Integrability of Lie algebroids} The
integrability of $A$ and the representability of $\cG(A)$ are not
exactly the same, due to the presence of isotropy groups. But,
since holonomy groupoids are always effective \cite{moerdijk}, we will
show that the
integrability of $A$ is equivalent to the representability of
$\cH(A)$.

\begin{pdef}[orbit spaces]
Let $\cX$ be a differentiable stack presented by a  Lie groupoid
$X=(X_1\rightrightarrows X_0)$. The orbit space of
$\cX$ is defined as the topological quotient $X_0/X_1$.
Throughout the paper, when we mention the orbit space is a
smooth manifold, we mean it has the natural smooth manifold
structure induced from $X_0$ (i.e. the projection $X_0\to X_0/X_1$ is
smooth).
\end{pdef}

\begin{proof}
We have to show the topological quotient is independent of choice
of presentations. Suppose that there is another presentation $Y$
which is Morita equivalent to $X$ through $(E, J_X, J_Y)$. Let
$O_x$ be the orbit of $X_1$ in $X_0$ through point $x$. By the
fact that both groupoid actions are free and transitive
fiber-wise, $J_Y\circ J_X^{-1}(O_x)$ is another orbit $O_y$ of
$Y$. In this way, there is a 1-1 correspondence between orbits of
$X$ and $Y$. Hence, $Y_0/Y_1$ understood as the space of orbits is
the same as $X_0/X_1$.
\end{proof}

\begin{thm} \label{integ-g} A Lie algebroid $A$ is integrable in the classical sense,
  i.e. there is a Lie groupoid whose Lie algebroid is $A$, iff the
  orbit space of $\cG(A)$ is a smooth manifold. Moreover, in this
  case the orbit space of $\cG(A)$ is the unique source-simply
  connected Lie groupoid integrating $A$.
\end{thm}
\begin{proof}First, let $Mon(P_a A)$ be the
monodromy groupoid of the foliation introduced by homotopy of
$A$-paths in Section \ref{apath}. We will show that $Mon(P_a A)$
is Morita equivalent to $Mon(P_0 A)$. Notice that $P_0A$ is a
submanifold of $P_a A$, so there is another groupoid $Mon(P_a
A)|_{P_o A}$ over $P_0 A$. We claim it is the same as $Mon(P_0
A)$. Namely, an A-homotopy $a(\epsilon, t)$ between two $A_0$
paths $a_0$ and $a_1$ can be homotopic to an $A_0$-homotopy
$\tilde{a} (\epsilon, t)$ between $a_0$ and $a_1$. The idea is to
divide $\tilde{a}$ into three parts:
\item{i)} First deform $a_0$ to $a_0^\tau$ through $a_0(\epsilon, t)$
  which is defined as $(1-\epsilon
  +\epsilon \tau' (t)) a_0((1-\epsilon) t +\epsilon \tau(t))$, where $\tau$ is the
  reparameterization induced in Section \ref{apath};
\item{ii)} Then, deform $a_0^\tau$ to $a_1^{\tau}$ through $a(\epsilon,
  t)^\tau$;
\item{iii)} Lastly, connect $a_1^{\tau}$ to $a_1$ through
  $a_1(\epsilon, t)$ which is defined as $a_1((1-\epsilon)
  \tau'(t) +\epsilon) a_1 (\epsilon t +(1-\epsilon) \tau (t))$.
  Then connect those three pieces by a similar method in the
  construction of concatenation (though it might be only piecewise smooth
  at the joints). Obviously, $\tilde{a}$ is a homotopy
  through $A_0$-paths and it is homotopic to a rescaling (over $\epsilon$) of $a(\epsilon, t)$ through
the concatenation of $a_0((1-\lambda)\epsilon, t)$,
  $(\lambda+(1-\lambda)\tau'(t))a(\epsilon,
  \lambda+(1-\lambda)\tau'(t))$ and $a_1((1-\lambda) \epsilon +
  \lambda, t)$. And eventually, we can smooth out everything to
  make the homotopy and the homotopy of homotopy both smooth so that
  they are as desired.

Then, it is routine to check that  $Mon(P_a A)|_{P_0 A}$ is Morita equivalent to $Mon(P_a A)$
through $\bt^{-1} (P_0 A)$, where $\bt$ is the target of the new
groupoid $Mon(P_aA)$.

So the orbit space of $\cG(A)$ can be realized as $P_a A/
Mon(P_a A)$. According to the main result in \cite{cf}, $P_a A/
Mon(P_a A)$ is a smooth manifold iff $A$ is integrable and if so,
$P_a A/ Mon(P_a A)$ is the unique source-simply connected Lie
groupoid integrating $A$.
\end{proof}

\begin{proof} [Proof for Theorem \ref{integ}]
First of all, by the same argument given in the proof above, one
can see that $Hol(P_0 A)=Hol(P_aA)|_{P_0A}$. Hence, $Hol(P_0 A)$is
Morita equivalent to $Hol(P_a A)$.

Moreover, if the orbit space of a holonomy groupoid  is a
manifold then it is Morita equivalent to the holonomy groupoid
itself (see \cite{moerdijk}).

Hence a differentiable stack $\cX=BG$  presented by a
holonomy groupoid $G$ is representable if and only if the orbit
space  $G_0/ G_1$ is a smooth manifold. One direction is obvious
because$G_0/G_1 \rightrightarrows G_0/G_1 $ is Morita equivalent
to $G=(G_1 \rightrightarrows G_0)$ if the orbit space is a
manifold. The converse direction is not hard to establish by examining the Morita equivalence diagram of $G$ and $\cX
\rightrightarrows \cX$. The Morita bibundle has to be $G_0$ since
$\cX$ is a manifold. Therefore $G_0$ is a principal $G$ bundle
over $\cX$. And this implies $G_0/G_1$ is the manifold $\cX$.

Notice that in general, the orbit spaces of monodromy groupoids
and holonomy groupoids of a foliation are the same. By Theorem
\ref{integ-g} and argument above, we conclude that  $A$ is
integrable iff $\cH(A)$ is representable and in this case,
$\cH(A)$ is $P_aA/Hol(P_a A)$, the unique source-simply connected
Lie groupoid integrating $A$.
\end{proof}

Combining the proofs of Theorem \ref{integ-g} and Theorem \ref{integ},
Theorem \ref{cf} follows naturally.

So far we have constructed $\cG(A)$ and $\cH(A)$ for every Lie
algebroid $A$ and verified that they are Weinstein groupoids.
Basically, we have done half of Theorem \ref{lieIII}. For the other
half of the proof, we first introduce some properties
of Weinstein groupoids. Before doing so, we give an example.

\begin{ep}[$B\Z_2$] $B\Z_2$ is a Weinstein group (i.e. its
base space is a point) integrating the trivial Lie algebra
0. The \'etale differentiable stack $B\Z_2$ is presented by
$\Z_2\rightrightarrows pt$ (here $pt$ represents a point).
We establish all the structure maps on this presentation.

The source and target maps are just projections from $B\Z_2$ to a
point. The multiplication $m$ is defined by
\[m: (\Z_2 \rightrightarrows pt)\times (\Z_2 \rightrightarrows pt)
\to (\Z_2\rightrightarrows pt), \; \text{by} \; m(a, b)= a\cdot
b,\] where $a, b\in \Z_2$. Since $\Z_2 $ is commutative, the
multiplication is a groupoid homomorphism (hence gives rise to a
stack homomorphism). It is easy to see that $m\circ (m\times id)=m
\circ (id \times m)$, i.e. we can choose the 2-morphism$\alpha$
inside the associativity diagram to be $id$.

The identity section $e$ is defined by
\[ e: (pt \rightrightarrows pt)\to (\Z_2 \rightrightarrows
pt), \; \text{by} \; e(1) =1, \] where 1 is the identity
element in the trivial group $pt$ and $\Z_2$.The inverse $i$ is defined by
\[i:(\Z_2 \rightrightarrows pt) \to (\Z_2
\rightrightarrows pt), \; \text{by} \; i(a)=a^{-1},\] where
$a\in \Z_2$. It is a groupoid homomorphism because $\Z_2$ is
commutative.

It is routine to check these maps satisfy the axioms of Weinstein
groupoids. The local Lie groupoid associated to $B\Z_2$ is just a
point. Therefore the Lie algebra of $B\Z_2$ is 0. Moreover, notice
that we have only used the commutativity of $\Z_2$. So for any
discrete commutative group $G$, $BG$ is a Weinstein group with Lie
algebra $0$.
\end{ep}

\begin{ep} [``$\Z_2 * B \Z_2$''] This is an example in which case Proposition
\ref{3asso} does not hold. Consider the groupoid $\Gamma=(\Z_2\times
\Z_2\rightrightarrows \Z_2)$. It is an action groupoid with
trivial $\Z_2$-action on $\Z_2$. We claim that the presented
\'etale differential stack $B\Gamma$ is a Weinstein group. We
establish all the structure maps on the presentation $\Gamma$.

The source and target maps are projections to a point.The
multiplication $m:\Gamma \times \Gamma \to \Gamma$ is defined by
\[m((g_1, a_1),(g_2, a_2))=(g_1 g_2, a_1 a_2).\]
It is a groupoid morphism because $\Z_2$ (the second copy) is
commutative. We have $m\circ(m\times id) = m\circ (id \times m) $.
But we can construct anon-trivial 2-morphism $\alpha:
\Gamma_0(=\Z_2) \times \Gamma_0 \times\Gamma_0 \to\Gamma_1$
defined by \[\alpha ( g_1, g_2, g_3) = ( g_1\cdot g_2\cdot g_3,
g_1\cdot g_2 \cdot g_3).\] Since the $\Z_2$ action on $\Z_2$ is
trivial, we have $m\circ (m\times id) = m\circ (id \times m)\cdot
\alpha $.

The identity section $e$ is defined by
\[ e: pt \rightrightarrows pt \to  \Gamma, \; \text{by}\;
e(pt)=(1, 1),\] where 1 is the identity element in $\Z_2$.The inverse
$i$ is defined by
\[ i: \Gamma \to \Gamma, \; \text{by} \; i(g, a)=(g^{-1},
a^{-1}). \] It is a groupoid morphism because $\Z_2$ (the second
copy) is commutative.

It is not hard to check $B\Gamma$ with these structures maps is a
Weinstein group. But when we look into the further obstruction of
the associativity described in Proposition \ref{3asso}, we found
failure. Let $F_i$'s be the six different ways of composing four
elements as defined in Proposition \ref{3asso}, then
the2-morphisms $\alpha_i$'s (basically coming from $\alpha$)
satisfy,
\[ F_{i+1} = F_i \cdot \alpha_i, \quad i=1, ..., 6 \;(F_7=F_1).\]
But $\alpha_i (1, 1, 1, -1)=(-1, -1)$ for all $i$'s except that
$\alpha_2 =id$. Therefore $\alpha_6 \circ \alpha_5\circ
...\circ\alpha_1(1,1,1,-1)=(-1, -1)$, which is not $id(1, 1, 1,
-1) = (-1, 1)$.
\end{ep}

\section{Weinstein groupoids and local groupoids}
In this section, we exam the relation between abstract
Weinstein groupoids and local groupoids.


Let us first show a useful lemma.
\begin{lemma}\label{eie} 
Given any \'etale atlas $G_0$ of $\cG$, there exists an open
covering $\{ M_l\}$ of $M$ such that the immersion $\bar{e}: M \to
\cG$ can be lifted to  embeddings $e_l: M_l \to G_0$. On the
overlap $M_l\cap M_j$, there exist an isomorphism $\varphi_{lj}$:
$e_j(M_j\cap M_l)\to e_l(M_j\cap M_l)$, such that $\varphi_{lj}
\circ e_j =e_l$ and $\varphi_{lj}$'s satisfy cocycle conditions.
\end{lemma}
\begin{proof}
Let $(E_e, J_M, J_G)$ be the HS-bibundle presenting the immersion
$\bar{e}: M \to \cG$. As  a right $G$-principal bundle over $M$,
$E_e$ is locally trivial, i.e. we can pick an open covering
$\{M_l\}$ so that $J_M$ has a section $\tau_l: M_l \to
E_e$ when restricted to $M_l$. Since $\bar{e}_l:=\bar{e}|_{M_l}$ is an immersion (the
composition of immersions $M_l\to M$ and $\bar{e}$ is still an
immersion), it is not hard to see that $pr_2: M_l \times_{\cG} G_0
\to G_0$ transformed by base change $G_0 \to \cG$ is an immersion.
Notice that $e_l=J_G \tau_l: M_l \to G_0$ fits inside a similar
diagram as \eqref{eq: immersion}:
\[
\begin{diagram}
\node{M_l\times_{\cG}G_0}\arrow{e,t}{pr_2}\arrow{s,t}{pr_1}
\node{G_0}\arrow{s,r}{\pi} \\
\node{M_l}\arrow{e,t}{\bar{e}_l}\arrow{ne,t}{e_l} \node{\cG}
\end{diagram}
\]
Following a similar argument in the proof of Lemma \ref{e-immersion}, we can find
a map $\alpha: M_l\times_{\cG}G\to G_1$ such that
\[ e_l\circ pr_1= pr_2 \cdot \alpha.\]
Since $\pi$ is \'etale, so is $pr_1$. Therefore $e_l$ is an
immersion.

Since an immersion is locally an embedding, we can choose an open
covering $M_{ik}$ of $\{ M_l\}$ so that $e_l|_{M_{ik}}$ is
actually an embedding. To simplify the notation, we can choose a
finer covering $\{ M_l\}$ at the beginning and make $e_l$ an
embedding. Moreover, using the fact that $G$ acts on $E_e$
transitively (fiberwise), it is not hard to find a local bisection
$g_{lj}$ of $G_1:=G_0\times_{\cG} G_0$, such that $e_l \cdot
g_{lj} =e_j$. Then $\varphi_{lj} = \cdot g^{-1}_{lj}$ satisfies
that $\varphi_{lj} \circ e_j = e_l$. Since $e_l$'s are embeddings,
$\phi_{lj}$'s naturally satisfy the cocycle condition.
\end{proof}

Before the proof of Theorem \ref{local}, we need a local statement.

\begin{thm}\label{local-local}
For every Weinstein groupoid $\cG$, there exists an open covering
$\{ M_l\}$ of $M$ such that one can associate a local Lie groupoid
$U_l$ over each open set $M_l$.
\end{thm}
\begin{proof} Let $\cG$ be presented by $G=(G_1\rightrightarrows G_0)$, and
$\{M_l\}$ be an open covering as in Lemma \ref{eie}.Let $(E_m, J_1, J_2)$
be the HS bibundle from $G_1\times_M G_1
\rightrightarrows G_0 \times_M G_0$ to $G$ which presents the
stack morphism $\bm: \cG\times_{M} \cG \to \cG$.
Notice that $M$ is the identity section, i.e.
\[
\begin{array}{ccc}
M_l \times_M M_l(=M_l) & \overset{ \bm=id}{\lra} &  M_l \\
\downarrow &\curvearrowright & \downarrow \\
\cG\times_M\cG &\overset{\bm} {\lra} &\cG.
\end{array}
\]
Translate this commutative diagram into groupoids. Then the
composition of HS morphisms
\begin{equation}\label{2hs}
\xymatrix{M_l\times_M M_l(=M_l)\ar[r]\ar[dd]\ar[dd]& G_1\times_M G_1\ar[dd]\ar@<-1ex>[dd]&  &G_1\ar[dd]\ar@<-1ex>[dd]\\
& & E_m\ar[dr]^{J_2}\ar[dl]_{J_1}& &\\
M_l\times_M M_l\ar[r]^{e_l\times e_l}&G_0\times_M G_0&  &G_0\\
}
\end{equation}
is the same (up to a 2-morphism) as $e_l:  M_l \to G_0$.
Therefore, composing the HS maps in \eqref{2hs} gives an HS
bibundle $J^{-1}_1 (e_l\times e_l (M_l\times_M M_l))$, which is
isomorphic (as an HS bibundle) to $M_l\times_{G_0} G_1$ which
represents the embedding $e_l$. Therefore, one can easily find a
global section
\[ \sigma_l: M_l\to  M_l\times_{G_0} G_1 \cong J_1^{-1}(e_l\times e_l
(M_l\times_M M_l))\subset E_m  \]  defined by $
x\mapsto (x, 1_{e_l(x)})$. Furthermore, we have $J_2 \circ
\sigma_l (M_l) =e_l(M_l)$. Since $G$ is an \'etale groupoid, $E_m$
is an \'etale principal bundle over $G_0\times_M G_0$. Hence
$J_1$ is a local diffeomorphism.
Therefore, one can choose two open neighborhoods $V_l\subset U_l$
of $M_l$ in $G_0$ such that there exists a unique section
$\sigma'_l$ extending $\sigma_l$ over $(M_l=M_l\times_M M_l
\subset) V_l\times_{M_l} V_l$ in $E_m$ and the image of $J_2\circ
\sigma'_l$ is $U_l$. The restriction of $\sigma'_l$ on $M_l$ is
exactly $\sigma_l$. Since $U_l\rightrightarrows U_l$ acts freely
and transitively fiberwise on $\sigma'_l(V_l\times_{M_l} V_l)$
from the right,  $\sigma'_l(V_l\times_{M_l} V_l)$ can serve as
an HS bibundle from $V_l\times_{M_l} V_l$ to $U_l$. (Here, we view
manifolds as groupoids.) In fact, it is the same as the morphism
\[m_l:=J_2\circ \sigma'_l: V_l \times_{M_l} V_l \to U_l.\]

 By a similar method, we can define the inverse as
follows. By (3), (4) and (5) in Definition \ref{wgpd}, we have
$\bar{i}\circ \bar{e}_l=\bar{e}_l$, so  the following diagram
commutes:
\[
\begin{array}{ccc}
M_l  & \overset{ \bm=id}{\lra} &  M_l \\
\downarrow &\curvearrowright & \downarrow \\
\cG &\overset{\bar{i}} {\lra} &\cG.
\end{array}
\]
Suppose $(E_i, J_1, J_2)$ is the HS bibundle representing
$\bar{i}$.  Translate the above diagram into groupoids, we have
the composition of the following HS morphisms:
\begin{equation} \label{2hs-i}
\xymatrix{M_l\ar[r]\ar[dd]\ar@<-1ex>[dd]& G_1\ar[dd]\ar@<-1ex>[dd]&  &G_1\ar[dd]\ar@<-1ex>[dd]\\
& & E_i\ar[dr]^{J_2}\ar[dl]_{J_1}& &\\
M_l\ar[r]^{e_l}&G_0&  &G_0\\
}
\end{equation}
is the same (up to a 2-morphism) as  $e_l:  M_l \to G_0$.
Therefore, composing the HS maps in \eqref{2hs-i} gives an HS
bibundle $J^{-1}_1 (e_l(M_l))$ which is isomorphic (as an HS
bibundle) to $M_l\times_{G_0} G_1$ which represents the embedding
$e_l$. Therefore, one can easily find a global section
\[ \tau_l: M_l\to  M_l\times_{G_0} G_1 \cong J_1^{-1}(e_l(M_l))\subset
E_i  \] defined by $
x\mapsto (x, 1_{e_l(x)})$. Furthermore, we have $J_2 \circ
\sigma_l (M_l) =e_l(M_l)$. Since $G$ is an \'etale groupoid, $E_i$
is an \'etale principal bundle over $G_0$. Hence $J_1$ is a local
diffeomorphism. Therefore, one can choose an open neighborhood of
$M_l$ in $G_0$, which we might assume as $U_l$ as well, such that there
exists a unique section $\tau'_l$ extending $\tau_l$ over $(M_l
\subset) U_l$ in $E_i$ and the image of $J_2\circ \tau'_l$ is in
$U_l$. The restriction of $\tau'_l$ on $M_l$ is exactly $\tau_l$.
So we can define
\[i_l:=J_2\circ \tau'_l: U_l \to U_l.\]

Since $M$ is a manifold, examining the groupoid picture of maps
$\bbs$ and $\bbt$, one finds that they actually come from two maps
$\bs$ and $\bt$ from $G_0$ to $M$. Hence, we define source and
target maps of $U_l$ as the restriction of $\bs$ and $\bt$ on
$U_l$ and denote them by $\bs_l$ and $\bt_l$ respectively.

The 2-associative diagram of $\bm$ tells us that $m_l\circ
(m_l\times id)$ and $m_l \circ (id\times m_l)$ differ in the
following way: there exists a smooth map from an open subset of
$V_l\times_{M_l} V_l \times_{M_l} V_l$, over which both of the above
maps are defined, to $G_1$ such that \[m_l\circ (m_l\times id)=
m_l \circ (id\times m_l) \cdot \alpha. \] Since the 2-morphism in
the associative diagram restricting to $M$ is $id$, we have
\[\alpha(x, x, x) =1_{e_l(x)}.\] Since $G$ is \'etale and $\alpha$ is
smooth, the image of $\alpha$ is inside the identity section of $G_1$.
Therefore $m_l$ is associative.

It is not hard to verify other groupoid properties in a similar
way by translating corresponding properties on $\cG$ to $U_l$.
Therefore, $U_l$ with maps defined above is a local Lie groupoid
over $M_l$.
\end{proof}

To prove the global result, we need the following Proposition:

\begin{prop}\label{glue}
Given $U_l$ and $U_j$ constructed as above (one can shrink
them if necessary), there exist an isomorphism of local Lie groupoids
$\tilde{\varphi}_{lj}: U_j \to U_l$ extending the isomorphism
$\varphi_{lj}$ in Lemma \ref{eie}. Moreover
$\tilde{\varphi}_{lj}$'s also satisfy cocycle conditions.
\end{prop}
\begin{proof}
Since we will restrict the discussion on $M_l\cap M_j$, we may
assume that $M_l=M_j$. According to Lemma \ref{eie},
there is a local bisection $g_{lj}$ of $G_1$ such that $e_l
\cdot g_{lj} =e_j$. Extend the bisection $g_{lj}$ to $U_l$ (we
denote the extension still by $g_{lj}$, and shrink $V_k$ and $U_k$
if necessary for $k=l, j$) so that \[ (V_l \times_{M_l} V_l) \cdot
(g_{lj} \times g_{lj} ) = V_j \times_{M_j} V_j \quad \text{and} \;
U_l \cdot g_{lj} = U_j.\] Notice that since $G_1$ is \'etale, the
source map is an local isomorphism. Therefore, by choosing small
enough neighborhoods of $M_l$'s, the extension of $g_{lj}$ is
unique. Let $\tilde{\varphi}_{lj}=\cdot g^{-1}_{lj}$. Then it is
naturally an extension of $\varphi_{lj}$. Moreover, by
uniqueness of the extension, $\tilde{\varphi}_{lj}$'s
satisfy cocycle conditions as $\varphi_{lj}$'s do.

Now we show $\tilde{\varphi}_{lj}=\cdot g_{lj}$ is a morphism
of local groupoids. It is not hard to see that $\cdot g_{lj}$
preserves source, target and identity embeddings. So we only have
to show that
\[ i_l \cdot g_{lj} =i_j,  \qquad m_l\cdot g_{lj}=m_j. \] For this purpose, we have
to recall the construction of these two maps. $i_l$ is defined as
$J_2 \circ \tau'_l$. Since there is a global section of $J_1$ over
$U_l$ in $E_i$, we have $J_1^{-1}(U_l)\cong U_l\times_{i_l, G_0}
G_1$ as $G$ torsors. Under this isomorphism, we can write
$\tau'_l$ as
\[ \tau'_l (x) = (x, 1_{e_l(x)}).\]
The $G$ action on $U_l\times_{i_l, G_0} G_1$ gives  $(x,
1_{e_l(x)})\cdot g_{lj}=(x, g_{lj})$. Moreover, we have
\[J_2((x, g_{lj}))=J_2(x, 1_{e_j(x)})=\bs_G(g_{lj}), \]
where $\bs_G$ is the source map of $G$. Combining all these, we
have shown that $i_l \cdot g_{lj} =i_j$. The other identity for
multiplications follows in a similar way.
\end{proof}

\begin{proof} [Proof of Theorem \ref{local}]
Now it is easy to construct $G_{loc}$ as in the statement of the
theorem. Notice that the set of $\{U_l\}$ with isomorphisms
$\varphi_{lj}$'s which satisfy cocycle conditions serve as a chart
system. Therefore, after gluing them together, we arrive at a
global object $G_{loc}$. Since $\varphi_{lj}$'s are isomorphisms
of local Lie groupoids, the local groupoid structures also glue
together. Therefore $G_{loc}$ is a local Lie groupoid.

If we choose two different open covering $\{ M_l\}$ and $\{M'_l\}$
of $M$ for the same \'etale atlas $G_0$ of $\cG$, we will arrive
at two systems of local groupoids $\{U_l\}$ and $\{U'_l\}$. Since
$\{M_l\}$ and $\{M'_l\}$ are compatible chart systems for $M$,
combining them and using Proposition \ref{glue}, $\{U_l\}$ and
$\{U'_l\}$ are compatible chart system as well. Therefore they
glue into the same global object up to isomorphisms near the
identity section.

If we choose two different \'etale atlases $G'_0$ and $G''_0$ of
$\cG$, we can take their refinement $G_0=G'_0\times_{\cG} G''_0$
and we can take a fine enough open covering $\{M_l\}$ so that it
embeds into all three atlases. Since $G_0\to G'_0$ is an
\'etale covering, we can choose $U_l$'s in $G'_0$ small enough so
that they still embed into $G_0$. So the groupoid constructed from the presentation $G_0$ with the
covering $U_l$ is the same as the groupoid constructed
from the presentation $G'_0$  with the covering $U_l$'s.   The same
is true for $G''_0$ and $G_0$. Therefore our local groupoid
$G_{loc}$ is canonical.

We will finish the proof of the Lie algebroid part in the next
section.
\end{proof}

\section{Weinstein groupoids and Lie algebroids}
In this section, we define the Lie algebroid of a Weinstein
groupoid $\cG$. An obvious choice is to define the Lie algebroid of$\cG$ as the Lie
algebroid of the local Lie groupoid $G_{loc}$.We give an
equivalent definition in a more direct way.

\begin{pdef} Given a Weinstein groupoid $\cG$ over $M$, there is
a canonically associated Lie algebroid $A$ over $M$.
\end{pdef}
\begin{proof}
We just have to examine more carefully the second part of proof
of Theorem \ref{local}. Choose an \'etale groupoid presentation
$G$ of $\cG$ and an open covering $M_l$'s as in Lemma \ref{eie}.
According to Theorem \ref{local}, we have a local groupoid $U_l$
and its Lie algebroid $A_l$ over each $M_l$. Differentiating the
$\tilde{\varphi}_{lj}$'s in Proposition \ref{glue}, we can achieve
algebroid isomorphisms $T\tilde{\varphi}_{lj}$'s which also
satisfy cocycle conditions. Therefore, using these data, we can
glue the $A_l$'s into a vector bundle $A$. Moreover, since the
$T\tilde{\varphi}_{lj}$'s are Lie algebroid isomorphisms, we can
also glue the Lie algebroid structures. Therefore $A$ is a Lie
algebroid.

Following the same arguments as in the proof of Theorem
\ref{local}, we can show uniqueness.  If we choose a different
presentation $G'$ and a different open covering $M_l$, we can just
choose the refinement of these two systems and will arrive at a
Lie algebroid which is glued from a refinement of both systems.
Therefore this is isomorphic to both  Lie algebroids constructed from
these two systems. Hence the construction is canonical.
\end{proof}

Now it is easy to see the following proposition holds:
\begin{prop}
Given a Weinstein groupoid $\cG$, it has the same Lie algebroid as
its associated local Lie groupoid $G_{loc}$.
\end{prop}

Together with the Weinstein groupoid $\cG(A)$ we have constructed
in Section \ref{w},  we are now ready to complete the proof of
Theorem \ref{lieIII}.

\begin{proof}[Proof of the second half of Theorem \ref{lieIII}]
We take the \'etale presentation $P$ of $\cG(A)$  and $\cH(A)$ as we
constructed in Section \ref{apath}. Let us recall how we construct
local groupoids from $\cG(A)$ and $\cH(A)$.

In our case, the HS morphism corresponding to $\bar{m}$ is $$(E:=
\bt_M^{-1}( m(P\times_M P)\cap \bs_M^{-1}(P), m^{-1}\circ \bt_M,
\bs_M).$$ The section $\sigma: M\to E$ is given by $x\mapsto
1_{0_x}$. Therefore if we choose two small enough open
neighborhoods $V\subset U$ of $M$ in $P$, the bibundle representing
the multiplication $m_V$ is a section $\sigma'$ over $V\times_M V$
of the map $m^{-1}\circ \bt_M$ in $E$.

Since the foliation $\cF$ intersects each transversal slice only
once, we can choose an open neighborhood $O$ of $M$ inside $P_0A$
so that the leaves of the restricted foliation $\cF|_O$ intersect
$U$ only once. We denote the homotopy induced by $\cF|_O$ as
$\sim_O$ and the holonomy induced by $\cF_O$ by $\sim_O^{hol}$.
Then there is a unique element $a \in U$ such that $a\sim_O
a_1\odot a_2$. Since the source map of $\Gamma$ is \'etale, there
exists a unique arrow $g: a_1\odot a_2 \curvearrowleft a$ in
$\Gamma$ near the identity arrows at $1_{0_x}$'s.

Then we can choose the section $\sigma'$ near $\sigma$ to be
\[\sigma': (a_1, a_2) \mapsto g.\]
So the multiplication $m_V$ on $U$ is
\[ m_V(a_1, a_2) = a(\sim_O a_1\odot a_2). \]

Because the leaves of $\cF$ intersect $U$ only once, $a$ has to be
the unique element in $U$ such that $a\sim^{hol}_O a_1\odot a_2$. It is not
hard to verify that both Weinstein groupoids give the same
local Lie groupoid structure on $U$.

Moreover, $U= O/\sim_O$ is exactly the local groupoid constructed
in Section 5 of \cite{cf}, which has Lie algebroid $A$. Therefore,
$\cG(A)$ and $\cH(A)$ have the same local Lie groupoid and their
Lie algebroids are both $A$.
\end{proof}

\begin{acknowledgements}We thank K. Behrend, H. Bursztyn,  M. Crainic, T. Graber,
 D. Metzler, I. Moerdijk, J. Mr{\v{c}}un, A. Weinstein, P. Xu and
 Marco Zambon for very helpful discussions and suggestions.
\end{acknowledgements}

\bibliographystyle{alpha}
\bibliography{bibz}

\end{document}